\RequirePackage{fix-cm}
\documentclass[smallextended]{svjour3}       % onecolumn (second format)
\smartqed  % flush right qed marks, e.g. at end of proof
\usepackage{graphicx}
\usepackage{algorithm}
\usepackage{algpseudocode}
\usepackage{amsmath}
\usepackage{enumerate}
\usepackage{subfigure}
\usepackage{amsmath}\DeclareMathOperator{\tr}{tr} \DeclareMathOperator{\diag}{diag}
\DeclareMathOperator{\Argmin}{Argmin}

\renewcommand\Im{\operatorname{Im}}
\journalname{Energy Systems}

\begin{document}

\title{Improved Spectral Clustering for Multi-Objective Controlled Islanding of Power Grid}
\titlerunning{Improved Spectral Clustering for Multi-Objective Controlled Islanding}        % if too long for running head

\author{Mikhail Goubko \and
        Vasily Ginz}


\institute{M. Goubko \at
              V.A. Trapeznikov Institute of Control Sciences of Russian Academy of Sciences, 117997, 65 Profsoyuznaya Street, Moscow, Russia, Tel.: +7-495-3349051\\
              \email{mgoubko@mail.ru}\\
           \and V. Ginz \at
                 Skoltech Center for Energy Systems, 143026, Skolkovo Innovation Center, 3 Nobel Street, Moscow, Russia, Tel.:
                 +7-495-2801481\\
              V.A. Trapeznikov Institute of Control Sciences of Russian Academy of Sciences, 117997, 65 Profsoyuznaya Street, Moscow, Russia}

\date{Received: date / Accepted: date}

\maketitle
\begin{abstract}
We propose a two-step algorithm for optimal controlled islanding that partitions a
power grid into islands of limited volume while optimizing several criteria: high
generator coherency inside islands, minimum power flow disruption due to teared
lines, and minimum load shedding. Several spectral clustering strategies are used in
the first step to lower the problem dimension (taking into account coherency and
disruption only), and CPLEX tools for the mixed-integer quadratic problem are
employed in the second step to choose a balanced partition of the aggregated grid
that minimizes a combination of coherency, disruption and load shedding. A greedy
heuristic efficiently limits search space by generating starting solution for exact
algorithm. Dimension of the second-step problem depends only on the desired number
of islands $K$ instead of the dimension of the original grid. The algorithm is
tested on standard systems with $118$, $2383$, and $9241$ nodes showing high quality
of partitions and competitive computation time.

\end{abstract}

\section{Introduction}
\label{intro} A power grid is a complex technical system; hence, it
is prone to technical disturbances. Even direct losses from
infrastructure damage and blackouts are enormous, so, much attention
is paid to assuring system sustainability under probable external
and internal shocks. Controlled islanding is the process of
splitting an interconnected power grid into smaller electrically
independent parts. It is used as a last-resort effort to cope with
many technical disturbances including undamped oscillations, voltage
collapse, cascading trips, etc.
\cite{ahmed2003scheme,sun2003ImbalancecoherenceOBDD}.

The rationale behind the process of controlled islanding is that a
smaller grid is easier to stabilize: islands have not be
synchronized, low frequency oscillations are less likely to occur in
a small grid, and so on. Also, the islanding operation can isolate
ill parts of the system from healthy ones, and the blackout will be
localized at ill islands if not avoided.

So, controlled islanding is a last chance to keep the grid alive (at
least, partially), which is a desirable effect. But there are also some side effects.

First of all, controlled islanding requires a series of complex
actions performed with high accuracy and coordination. Any error can
cause a cascading trip of generators or transmission lines.

Secondly, some power transmission lines are switched off during the islanding
operation, which makes a great shock to the grid even when islanding is accurately
planned and perfectly implemented. Such disruption leaves a partitioned grid in a
highly unstable state, making it questionable to stabilize the state of some
islands. The simplest metrics of the power flow disruption is the total volume of
power flows broken during an islanding operation.

Thirdly, a well-designed interconnected grid has more opportunities
to serve the current demand than any partitioned grid due to limited
power transmission and ramp rate opportunities of the latter.
Therefore, some load shedding is an essential part of the controlled
islanding process
\cite{junchen2015CoherencyDisruptionImbalanceCompare}.

After all, restoring the grid after the controlled islanding is also
a time-consuming and complex operation.

Bearing in mind high risks of controlled islanding, it is important
to plan properly the islanding operation and to suggest a realistic
and safe grid partitioning scheme that will stabilize the grid and
minimize side effects.

Optimal islanding of a real-world power system is a high-dimensional
optimization problem. The islanding decision is made under the
extremal time pressure; a system operator has just several
seconds to develop and implement an islanding scheme. Therefore,
finding an optimal grid partition becomes a non-trivial task.

Metrics of system stability considered by power system security
studies are computationally expensive. Shocks incurred by islanding
operation must be classified as large-scale, and non-linear effects
cannot be neglected. Hence, to make reliable predictions of the
after-islanding system dynamics extensive time-domain simulations
are performed under the paradigm of transient stability analysis
\cite{kundur2004StabilitySurvey}.

\sloppy That is why most formal models of optimal controlled islanding (OCI) do not
consider directly restoring system stability as an optimization criterion.
Workarounds include incorporating external constraints (e.g., limited island volume)
or constructing some simpler criteria, e.g., degree of generators' coherency (or
dynamic coupling) in the pre-islanding system state. Many articles leave stability
constraints behind the scene (e.g.,
\cite{ding2013CoherencyDisruptionConstrainedSpectralBisection,quiros2014determination,sanchez2014AdmittanceDisruptionHierSpec,trodden2014optimization}
and many others) concentrating on minimization of side effects of the islanding
process, and most of them take into account only a single aspect (load shedding in
\cite{fanpardalos2012SLDC_MILP,pahwa2013ShedLoadDCGreedy}) or a couple of aspects
(e.g., generator coherency and flow disruption in
\cite{ding2013CoherencyDisruptionConstrainedSpectralBisection}).

A mathematical model of OCI introduced in the present article takes into account
multiple aspects of the islanding process and their corresponding performance
metrics (generator coherency \cite{chow1987Sparsetimescale}, flow disruption
\cite{li2005strategic}, load shedding, and, optionally, line susceptance). A fast
and efficient two-step algorithm is proposed to calculate a rational scheme of
partitioning a grid into the desired number of islands $K$ with the limited maximum
island volume.

In the first step hierarchical spectral clustering algorithms from
\cite{sanchez2014AdmittanceDisruptionHierSpec,quiros2015constrained}
are used to break the grid into $n'>K$ islands to minimize the
weighted combination of coherency and flow disruption metrics.

In the second step each of detailed grid partitions obtained in the
first step is transformed into the aggregated grid with $n'$
vertices. The aggregated grid is then partitioned into $K$ connected
islands with the limited maximum island volume by an exact algorithm
implemented in CPLEX 12. A greedy heuristics provides an efficient
starting solution, which sufficiently fosters calculations. In
addition to generator coherence and flow disruption the optimization
criterion in the second step also takes into account the power imbalance
inside islands.

The idea is to combine high speed (as the problem dimension is efficiently reduced
in the first step) and high flexibility (complex optimization criteria and
constraints are allowed in the second step). Computational experiments on the models
of standard power systems with 118, 2383, and 9241 nodes show that the proposed
algorithm is fast enough and outperforms alternative approaches both in bulk and in
the value of every single performance metric.

The rest of the article has the following structure. In Section \ref{sec:Review}
recent approaches to OCI are surveyed. Then in Section \ref{sec:Defs} we introduce
the notation and basic mathematical concepts used to define controlled islanding
performance metrics. In Section \ref{sec:Spectral} essential information is provided
about spectral clustering, which is the main tool for fast and efficient islanding.
The two-step algorithm for multi-objective grid partitioning is introduced in
Section \ref{sec:Algorithm}, while computational experiments on three cases of power
grids are presented in Section \ref{sec:Experiments}. Section \ref{sec:Conclusion}
concludes with some open issues and perspectives.

\section{Literature Review}\label{sec:Review}

After a critical breakdown (e.g., a shortage, a circuit trip, or a
generator failure) a power grid may become instable. Different
groups of generator go out of sync, and the main goal of the
automated grid control is to return the system into the stable state
with minimum load shedding. In
\cite{ahmed2003scheme,vittal1998FirstStabilityManual} it is shown that
controlled islanding (accompanied with appropriate load shedding)
can be a promising strategy to prevent cascading blackouts in power
systems. A sort of modal analysis (analysis of normal forms) was
employed in \cite{vittal1998FirstStabilityManual} to perform
generators' grouping while in \cite{ahmed2003scheme} predefined
islands were considered. At the same time, the choice of a rational
grid partitioning strategy was a challenging discrete optimization
problem. System stability and the amount of load shedding were
considered as the main criteria of partition quality.

\emph{Shed load} is the amount of load that cannot be served safely
given the topology of islands in the grid according to voltage and
safety constraints and, thus, should be disconnected during the
islanding operation. Shed load can be obtained by solving the
optimal load shedding (OLS) problem under the alternating current
(AC) model. AC-OLS reduces to a non-linear optimization problem
\cite{pahwa2013ShedLoadDCGreedy}, which takes sufficient time to
solve. A simpler version of OLS problem often used in contingency
analysis is the direct current (DC) approximation, which reduces to
the linear program with (real) power flow balance, phase angle and
maximum real flow constraints \cite{pahwa2013ShedLoadDCGreedy}.

Even in a small power grid the number of alternative islanding
schemes is enormous. An operator has just few seconds to choose and
implement a grid islanding scheme (e.g., some generators go out of
step within 5 sec after the contingency in the scenario modeled in
\cite{xu2010SimplifyMETIS}). A detailed time-domain simulation
cannot be run for every alternative partition to verify its
stability, and indirect stability indicators are typically used at
the partition selection stage (a notable exception is the use of
PSSENG time-domain simulation in the genetic algorithm
\cite{elWerfelli2008StabilityGenetic} to evaluate island stability
of IEEE 118-bus scheme).

A popular stability indicator of an electrical power system is coherency of its
generators \cite{you2004CoherencySimulation} (their aspiration to swing together).
Generator grouping methods based on slow coherency detection were primarily
developed for system model reduction
\cite{chow1987Sparsetimescale,chow1988Twotimescale,chow1991Twotimescale}. Several
approaches were proposed (see \cite{chow1995CoherencyGreedy,yusof1993Coherency}) to
construct complete grid partitions by assigning loads to generator groups, but they
ignored load shedding and the other metrics relevant to controlled islanding (e.g.,
transmission line security constraints).

In
~\cite{sun2003ImbalancecoherenceOBDD,zhao2003ImbalancecoherenceOBDDpostprocessing}
controlled islanding is considered as a satisfiability problem (the
problem of finding a grid partition that fulfills a set of
constraints). Ordered binary decision diagrams (OBDD) were used for
limited enumeration of islanding schemes that satisfy generator
coherency and load/generation balance constraints, while DC-OLS is
run for each candidate solution until stability constraints are
verified. OBDD-based enumeration is time-demanding, so large
real-world networks have to be aggregated before applying the
algorithm.

Load shedding is minimized in \cite{pahwa2013ShedLoadDCGreedy} by solving a series
of DC-OLS problems inside a greedy algorithm. Stochastic programming is used in
\cite{golari2014stochastic,golari2016large} to find an islanding scheme, which
minimizes the average load shedding against a series of pre-defined contingencies.
In
\cite{fanpardalos2012SLDC_MILP,trodden2013milp,trodden2014optimization,ding2015SLMILP}
the grid partitioning problem is reduced to the mixed-integer linear program (MILP)
with additional generator coherency constraints in
\cite{trodden2013milp,trodden2014optimization,ding2015SLMILP} and connectivity
constraints in \cite{ding2015SLMILP}. MILP is solved with exact algorithms
implemented in CPLEX numeric optimization package. Nevertheless, computational
experiments on grids with $\le 300$ nodes show that these algorithms are time
expensive and, hence, hardly applicable to online islanding problems.

To improve time efficiency of algorithms the amount of load shedding
is often approximated by the load/generation imbalance. In
particular, a heuristic algorithm is proposed in
\cite{wang2004ImbalanceCoherenceHeuristic} to minimize total
load/generation imbalance under generator coherency constraints, but
its time and cost efficiency was not compared to alternative
approaches.

In addition to generator coherency, the power flow disruption is
another popular indicator of the island (in)stability
\cite{li2005strategic,peiravi2009fast}. When a transmission line is
tripped during an islanding operation, the power flow through this
line immediately drops to zero. The greater is the disruption (total
power flow through the lines being teared), the greater is the
excessive shock to the islanded power system (in addition to the
initial disturbance) and the less possible is its stabilization.

Algorithms of controlled islanding based on spectral clustering
techniques are intensively studied in recent years. The idea of
using spectral clustering techniques for fast partitioning of
electrical networks can be traced back to
\cite{li2005strategic,yang2007UnclearSpectral}. A Spectral
Clustering Controlled Islanding (SCCI) scheme is proposed in
\cite{ding2013CoherencyDisruptionConstrainedSpectralBisection}. It
accounts both for generator coherency and power flow disruption. In
the first step the normalized spectral bisection is used to divide
the generators in two coherent groups. In the second step loads are
assigned to generator groups using the unnormalized constrained
spectral bisection to minimize the flow disruption. Recursive
bisection is applied to obtain the desired number of islands. SCCI
algorithm is shown to be much faster than OBDD-based enumeration
techniques
\cite{zhao2003ImbalancecoherenceOBDDpostprocessing,sun2003ImbalancecoherenceOBDD}
with the minor loss in solution quality.

Later this methodology was extended in different aspects. In
particular, $k$-medoids algorithm is used to cluster eigenvector
points \cite{ding2014constrained} instead of $k$-means suggested in
\cite{li2005strategic}. The problem of outliers when performing
eigenvector analysis is addressed in \cite{ding2014constrained}.
Recursive bisection were replaced by direct $k$-way partitions in
\cite{ding2014constrained,sanchez2014AdmittanceDisruptionHierSpec},
which decreased computation cost and increased partition quality.

In \cite{sanchez2014AdmittanceDisruptionHierSpec} it is proposed to
incorporate hierarchical clustering algorithm (first introduced in
\cite{ward1963hierarchical}) into the spectral clustering scheme to
stimulate generation of connected islands.
%Distance matrix for
%hierarchical clustering is built by embedding grid vertices onto the
%$k$-dimensional sphere using the coordinates from $k$ first
%eigenvectors of the Laplace matrix and using the shortest graph
%paths.
In \cite{quiros2015constrained} this approach is extended to
account for generator coherency. Predefined generator groups are
considered and loads are assigned to the nearest neighboring
generator in the same spectral graph embedding. Finally, in
\cite{quiros2014determination} the spectral clustering approach is
applied to the problem of parallel system restoration, which differs
from controlled islanding in several important aspects.

An essential limitation of OCI algorithms based on spectral
clustering is their disregard of load shedding. A yet another
serious gap is that spectral methods sometimes result in extremely
\emph{imbalanced} partitions (those consisting of a huge mainland
and several tiny islets), which is not practical in many ways.

The present article is devoted to further development of the spectral approach to
OCI. We overcome existing shortcomings by using the hierarchical spectral clustering
algorithms from \cite{sanchez2014AdmittanceDisruptionHierSpec,quiros2015constrained}
as pre-processing routines to decrease the problem dimensionality to the desired
degree caring only for generator coherency and flow disruption. Island volumes and
shed load are accounted for in the second step of the proposed algorithm, where
mixed-integer quadratic problem (MIQP) is solved.

\section{Basic Notation}\label{sec:Defs}
Every islanding decision is unique, since it is made under unique
(and often unexpected) conditions. The state of the power system at
the moment before the islanding decision is described by several
groups of variables. Table \ref{tab:nomenclature} summarizes the notation used to define the
context of an islanding process and its performance metrics.

Let us consider a power grid consisting of $n$ nodes (so called, \emph{buses})
indexed from 1 to $n$, and $m$ transmission lines. Without loss of generality assume
that at most one generator or a load is assigned to each bus and all generators are
located at the first $n_g\le n$ buses. \textbf{Bold} is used for vectors.

\begin{table}[!ht]
\caption{Nomenclature}
\label{tab:nomenclature}       % Give a unique label
\begin{tabular}{p{2.5cm} p{8.5cm} }
\hline\noalign{\smallskip}
Notation & Description \\
\noalign{\smallskip}\hline\noalign{\smallskip}
    $N=\{1,...,n\}$ & Set of buses in the grid\\
    $N_g=\{1,...,n_g\}$ & Set of buses with generators installed\\
    $V_i$ & Voltage amplitude at bus $i\in N$\\
    $\theta_i$ & Phase angle at bus $i\in N$\\
    $G_i$ & Maximum real power output of generator at bus $i\in N_g$\\
    $g_i$ & Current real power output of generator at bus $i\in N_g$\\
    $D_i$ & Real power demand at bus $i\in N$\\
    $d_i$ & Current real load at bus $i\in N$\\
    $p_i=d_i-g_i$ & Current real power injection at bus $i\in N$\\
    $p(s)=\sum_{i\in s} p_i$ & Load/generation imbalance in island $s\subseteq N$\\
    $p_{ij}$ & Current real power flow from bus $i\in N$ to bus $j\in N$, $i>j$, \\
             & ($p_{ij}<0$ if the flow is directed from $j$ to $i$) \\
    $P=(p_{ij})_{i,j=1}^n$ & Real power flow matrix (lower-triangular)\\
    $\bar{p}_{ij}$ & Real power flow limit for the line between buses $i,j\in N$\\
    $H_i$ & Inertia constant of machine at bus $i\in N_g$\\
    $Y=(y_{ij})_{i,j=1}^n$ & Nodal complex admittance matrix\\
    $B=(b_{ij})_{i,j=1}^n$ & Nodal susceptance matrix, $B=\Im(Y)=\left(%
\begin{array}{cc}
  B_{11} & B_{12} \\
  B_{21} & B_{22} \\
\end{array}%
\right)$, \\
& where $B_{11}$ is $n_g\times n_g$ sub-matrix limited to generator buses only\\
    $\Delta=(\Delta_{ij})_{i,j=1}^n$ & Electrical distance matrix defined in \cite{cotilla2013multiGA} on the basis of $B$ matrix\\
    $\tilde{B}=(\tilde{b}_{ij})_{i,j=1}^{n_g}$ & Reduced susceptance matrix \cite{butyrin2003}: $\tilde{B}=B_{11} - B_{12}B_{22}^{-1}B_{21}$\\
    $\tilde{\Phi}=(\phi_{ij})_{i,j=1}^{n_g}$ &  Dynamic coupling matrix of generators \cite{ding2013CoherencyDisruptionConstrainedSpectralBisection}, \\
    & where $\phi_{ij}=\left(\frac{1}{H_i}+\frac{1}{H_j}\right)|V_i||V_j|\tilde{b}_{ij}\cos
(\theta_i-\theta_j)$ \\
    $\Phi = \left(%
\begin{array}{cc}
  \tilde{\Phi} & 0 \\
  0 & 0 \\
\end{array}%
\right)$ & Dynamic coupling matrix ($n\times n$) of generator buses\\
    $\mathbf{w}=(w_i)_{i=1}^n$ & The vector of bus \emph{volumes}, where $\mathbf{w}=|P|\cdot \mathbf{1}_n$   (alternatively, $w_i=G_i+D_i$; another volume metric is also possible)\\
    $w(s)=\sum_{i\in s} w_i$ & Volume of island $s\subseteq N$\\
\noalign{\smallskip}\hline\noalign{\smallskip}
\end{tabular}
\end{table}

Let $\mathbf{1}_n$ be the $n$-dimensional all-ones vector,
$\mathbf{0}_n$ be the $n$-dimensional all-zeros vector, and let
$I_{n\times n}=\diag(\mathbf{1}_n)$ be the $n\times n$ identity
matrix. For real symmetric $n\times n$ matrix $A=(a_{ij})_{i,j=1}^n$
its \emph{Laplace matrix} is defined as
$$L(A):=\diag(A\cdot\mathbf{1}_n)-A.$$
    of \emph{volumes} $\mathbf{w}$ is
also given (see Table \ref{tab:nomenclature} for possible definitions of $\mathbf{w}$), the \emph{symmetric normalized Laplace matrix} of matrix
$A$ under volumes $\mathbf{w}$ is defined as
$$L_{sym}(A|\mathbf{w}):=\diag(\mathbf{w})^{-\frac{1}{2}}\cdot L(A)
\cdot \diag(\mathbf{w})^{-\frac{1}{2}}.$$

Power flows in a power grid are naturally modeled by a directed graph with weights
assigned to its vertices (representing injections in buses) and to arcs
(representing real power flows through transmission lines). Let us introduce the
notion of the \emph{graph cut}, which plays the central role throughout the article.
Consider a simple directed graph with vertex set $N$ and arc weights $a_{ij}$, $i,j
\in N$. For vertex set $s\subseteq N$ the \emph{cutset} is the minimum set of arcs
needed to be removed to isolate $s$ from the rest of the graph. The total weight of
the cutset is called the \emph{cut}, which is calculated as
$$Cut_A(s):=\sum_{i\in s,j\in N\backslash s}a_{ij}=x^T L(A) x,$$
where $A:=(a_{ij})_{i,j=1}^n$ is the matrix of arc weights, and $x$ is an indicatory
vector of vertex set $s$ (i.e., $x_i=1$ if vertex $i\in s$, and $x_i=0$ otherwise).

An \emph{islanding scheme} with $K$ islands is represented by
partition $\pi=(s_1,...,s_K)$ of vertex set $N$ into $K$ disjoint
parts. For partition $\pi=(s_1,...,s_K)$ the \emph{cut} is defined
as
$$Cut_A(\pi):=\sum_{k=1}^K Cut_A(s_k)=
\tr X^T L(A)X,$$ where $X$ is $n\times K$ indicatory matrix of
partition $\pi$ (i.e.,  $x_{ik}$ is equal to unity when $i\in s_k$
and is zero otherwise).

For positive vector $\mathbf{w}=(w_i)_{i=1}^n$ of vertex volumes the
\emph{normalized cut} is defined as
\begin{equation}\label{eq:ncut}
NCut_A(\pi|\mathbf{w}):=\sum_{k=1}^K \frac{Cut_A(s_k)}{w(s_k)}= \tr
Z^T L_{sym}(A|\mathbf{w})Z,
\end{equation}
where $w(s):=\sum_{i\in s} w_i$ is the volume of island $s\subseteq N$, $Z =
\diag(\mathbf{w})^{\frac{1}{2}} X \diag(X^T\mathbf{w})^{-\frac{1}{2}}$ (in the other
words, $z_{ik} = \sqrt{\frac{w_i}{w(s_k)}}$ if $i \in s_k$ and is zero otherwise.)
Matrix $Z$ is \emph{orthogonal}, i.e., $Z^TZ=I$.

Many popular performance metrics of controlled islanding can be
written in terms of a graph cut with appropriately weighted graph
arcs. For any $K$-partition $\pi=(s_1,...,s_K)$ \emph{generators'
dynamic coupling} is defined in
\cite{ding2013CoherencyDisruptionConstrainedSpectralBisection} as
\begin{equation}\label{eq:C_def}
C(\pi)=Cut_\Phi(\pi), \end{equation} and \emph{power flow
disruption} is defined as
\begin{equation}\label{eq:D_def}
D(\pi)=Cut_{|P|}(\pi).
\end{equation}
Here the \emph{absolute power flow matrix} (nonnegative and symmetric) is defined as
$|P|:=(\frac{1}{2}|p_{ij}|+\frac{1}{2}|p_{ji}|)_{i,j=1}^n$, where $p_{ij}$ is the
\emph{real power flow} along arc $ij$.

In the present article we follow
\cite{ding2013CoherencyDisruptionConstrainedSpectralBisection} and use \emph{dynamic
coupling} matrix $\Phi$ (defined in Table \ref{tab:nomenclature}) to build metric
(\ref{eq:C_def}) of generator coherence. Alternatively, more sophisticated
approaches (e.g., independent component analysis \cite{ariff2013coherency} or
hierarchical trajectory cluster analysis \cite{juarez2011characterization}) can be
employed to build another matrix of generators' coherence, which can be used in
$C(\pi)$.

The \emph{electrical cohesiveness index} (ECI), defined as
$$ECI(\pi):=Cut_\Delta(\pi),$$
was used in \cite{cotilla2013multiGA} to split the grid on the basis of electrical
distance matrix $\Delta$.

The \emph{shed load} is the amount of load that cannot be served
safely given the topology of islands in the grid according to
voltage and safety constraints. For partition $\pi=(s_1, ..., s_K)$
the shed load is denoted as
$$S(\pi):=\sum_{k=1}^K S(s_k),$$ where $S(s_k)$ is the load shed in
island $s_k$. Several approaches with different accuracy and
computational complexity are used to evaluate $S(s_k)$. Let
$S_{AC}(s_k)$ be the amount of load shedding in island $s_k$ in the
solution of AC-OLS. Since AC-OLS is computationally expensive, below
we calculate it only for final partitions to verify the quality of
the solution.

DC-OLS is a simpler version of OLS problem limited to real flows
with fixed bus voltages, zero line losses, and low phase angle
differences. DC-OLS reduces to the linear program with real power
flow balance, phase angle and maximum real flow constraints
\cite{pahwa2013ShedLoadDCGreedy}. The load shed in island $s_k$
calculated from OLS-DC solution is denoted by $S_{DC}(s_k)$.

The estimate of shed load calculated as a solution of the maximum
flow problem for island $s_k$ obtained by relaxing phase angle
constraints in DC-OLS, is denoted as $S_{MF}(s_k)$ (see more details
in Section \ref{sec:MIQP} below). Finally, when maximum real flow
constraints are relaxed, the minimum amount of shed load is
estimated by \emph{excess load}:
\begin{equation}\label{eq:I_def}
S_{EL}(s_k):=\max\left[p(s_k); 0\right],
\end{equation}
where $p(s):=\sum_{i\in s} (d_i-g_i)$ is total imbalance between load $d_i$ and
generation $g_i$ in nodes $i\in s$ of island $s\subseteq N$. If losses in lines are
neglected, then $p(s_k)=Cut_P(s_k)$, where $P=(p_{ij})_{i,j=1}^n$ is a matrix of
real power flows, and excess load is also expressed using the (directed) graph cut.

Finally, we note that  $S_{EL}(s) \le S_{MF}(s) \le S_{DC}(s) \le
S_{AC}(s)$ for any island $s\subseteq N$, so AC-OLS gives the most
conservative estimate of the shed load.

\section{Spectral Clustering Basics}\label{sec:Spectral}

Spectral clustering is an approach to finding approximate solutions
of minimum cut problems for weighted undirected graphs. Consider an
undirected simple graph $\langle N, E\rangle$ with vertex set
$N=\{1,...,n\}$, edge set $E\subseteq N\times N$, and non-negative
edge weights $a_{ij}$, $ij\in E$. Then, given matrix
$A=(a_{ij})_{i,j=1}^n$ (non-negative and symmetric), the graph
minimum $K$-cut (or $K$-partition) problem is to find a partition
$\pi=(s_1,...,s_K)$ of a vertex set $N$ into $K$ disjoint parts that
minimizes $Cut_A(\pi)$.\footnote{Alternatively, the minimum $K$-cut
problem is to minimize the total weight of edges, which, if removed,
break the graph into $K$ connected components. These two definitions
are, in fact, equivalent.}

If, in addition, positive volume $w_i$ is assigned to every vertex
$i\in N$, the minimum balanced $K$-cut problem is to minimize
$Cut_A(\pi)$ by choosing a partition $\pi=(s_1,...,s_K)$, such that
$w(s_k)\le W, k=1,...,K$, where $W \ge W(N)/K$ is an upper cluster
volume limit.

The minimum $K$-cut problem is known
\cite{goldschmidt1988polynomial} to be solvable in polynomial time
$O(n^{K^2})$ for any fixed $K$, but is NP-complete
\cite{GareyJohnson1979} when $K$ is not limited. The minimum
balanced cut problem is NP-complete \cite{GareyJohnson1979} even for
the graph bipartition problem (when $K=2$, $n$ is even, $w_i=1$, and
$W = n/2$).

Many efficient approximate graph partitioning algorithms were
developed since then, with spectral clustering being among most
popular ones. For spectral clustering the balanced $K$-cut problem
is replaced with $NCut_A(\pi|\mathbf{w})$ minimization problem with
no cluster volume constraints resulting is some sort of relaxation:
for partition $\pi = (s_1, ..., s_K)$
\begin{equation}\label{eq:NCut}
NCut_A(\pi|\mathbf{w}) =
\frac{Cut_A(s_1)}{w(s_1)}+...+\frac{Cut_A(s_K)}{w(s_K)},
\end{equation}
and denominators in (\ref{eq:NCut}) penalize small cluster volumes,
while for equal cluster volumes (such that $w(s_1)=...=w(s_K)$) we
have
$$NCut_A(\pi|\mathbf{w})= Cut_A(\pi)\cdot \frac{K}{w(N)}.$$

Spectral clustering is based on the spectral lower bound for the trace minimization
problem \cite{lutkepohl1996handbook}. For an arbitrary real symmetric $n\times n$
matrix $A$ let $\lambda_i(A)$, $i=1,...,n$, be its eigenvalues enumerated in
ascending order, and let $\mathbf{u}_i(A)$ be corresponding eigenvectors. Then for
any $n\times K$ orthogonal matrix $Z$, $K\le n$, the following inequality holds:
\begin{equation}\label{ineq:spectral}
\tr Z^T A Z \ge \lambda_1(A) + ... + \lambda_K(A) = \tr U^T A U,
\end{equation}
where $U = (\mathbf{u}_1(A), ..., \mathbf{u}_K(A))$ is a matrix
composed of $K$ first (normalized) eigenvectors of matrix $A$.

Any $K$-partition can be written in the form of an orthogonal
$n\times K$  binary matrix $X = (x_{ik})$, where $x_{ik} = 1$ if an
only if $i$-th node belongs to $k$-th cluster. Spectral clustering
algorithms approximate expression $\tr U^T A U$ in
(\ref{ineq:spectral}) with some admissible partition $X$, which is
close in some sense to matrix $U$. The convenient notation is
developed below for basic algorithms of spectral clustering, which
is used in the next section to build combined partitioning
algorithms.

According to equation (\ref{eq:ncut}), the normalized cut is written
with normalized Laplace matrix $L_{sym}(A|\mathbf{w})$, so in the
classical algorithm \cite{ng2002spectral} the rows of matrix $U$,
whose columns are the first $K$ eigenvectors of $L_{sym}(A|\mathbf{w})$,
are normalized to the unit Euclidian norm and then partitioned with
$k$-means clustering. Each row corresponds to a graph vertex. Denote
the resulting $K$-partition with $\kappa_K(A|\mathbf{w})$.
%In unnormalized spectral clustering the same procedure is applied to
%matrix $L(A)$ instead of $L_{sym}(A|\mathbf{w})$. Let $\kappa_K(A)$
%stand for the corresponding $K$-partition.

The algorithm of $k$-means is known to be sensitive to outliers
(arising, for example, when graph has pendent vertices). It is shown
in \cite{demetriou2013implementing} that more stable islanding
schemes for real power grids can be obtained when $k$-means is
replaced with the more robust $k$-medoids algorithm. The
corresponding partition is denoted with $\mu_K(A|\mathbf{w})$.
%$(\mu_K(A)$ for unnormalized spectral clustering).

Both $k$-means and $k$-medoids often suggest disconnected subgraphs as partition
elements. To deal with with problem it is suggested in
\cite{sanchez2014AdmittanceDisruptionHierSpec} to use hierarchical clustering
instead. In the Hierarchical Spectral Clustering (HSC) algorithm
\cite{sanchez2014AdmittanceDisruptionHierSpec} the graph is pre-processed by
iteratively merging pendent vertices to their neighbors. No pendent vertex are left
in the graph to avoid the outlier problem. Then eigenvector matrix $U$ is calculated
for the normalized Laplacian of the absolute power flow matrix $|P|$. Normalized rows of
matrix $U$ are considered as points on a $K$-dimensional sphere and the graph being
partitioned is embedded onto this sphere. The distance between any pair of graph
vertices is calculated as the length of the shortest path in the embedded graph (the
cosine distance between incident nodes is considered). A hierarchical clustering
algorithm \cite{ward1963hierarchical} is then applied to the obtained distance
matrix (Standard Matlab implementation with complete linkage
\cite{defays1977efficient} is used below.) The shortcoming of this algorithm is
bigger variation of island volumes: a partition often consists of one or two big
islands and a collection of small islets. Let $\chi_K(|P| ~ |\mathbf{w})$ stand for the
$K$-partition built with hierarchical normalized spectral clustering.
% and $\chi_K(A)$ -- for the corresponding unnormalized clustering partition.

\begin{figure}
\includegraphics[width=0.55\textwidth]{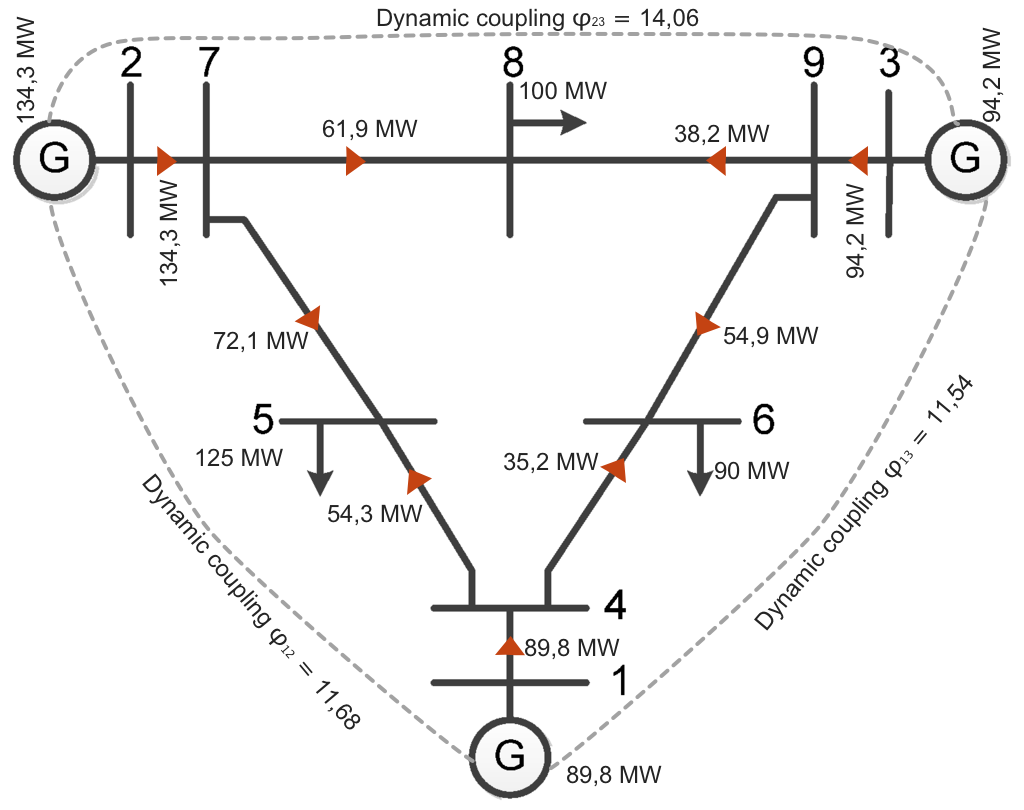}
\caption{Example of pre-islanded network and power flows}\label{fig:pre9bus}
\end{figure}

Let us consider a simplistic numeric example by partitioning a tiny power grid (the
IEEE 9-bus system shown in Figure \ref{fig:pre9bus}) by HSC algorithm. The network
has 3 generators installed, 3 loads and the total of 9 buses. The (outbound) real
power flows, current real power generations and loads calculated using AC-OPF are
presented in Figure \ref{fig:pre9bus} with arrows showing the flow direction. To
calculate a bisection of this network with HSC algorithm, pendent nodes 1,2, and 3
are joined respectively to nodes 4,7, and 9. The matrix of absolute power flows for
the resulting 6-nodal graph is
$$|P|=\left(%
\begin{array}{cccccc}
  0    & 27.15& 17.60& 0    & 0    & 0 \\
  27.15& 0    & 0    & 36.05& 0    & 0 \\
  17.60& 0    & 0    & 0    & 0    & 27.45 \\
  0    & 36.05& 0    & 0    & 30.95& 0 \\
  0    & 0    & 0    & 30.95& 0    & 19.10 \\
  0    & 0    & 27.45& 0    & 19.10& 0 \\
\end{array}%
\right),$$ and matrix
$$U=\left(%
\begin{array}{cc}
-0.4603 &  -0.0075\\
-0.2665 &   0.3019\\
-0.4708 &   0.7119\\
-0.2808 &  -0.0240\\
-0.5631 &  -0.5882\\
-0.3155 &  -0.2355\end{array}\right)$$ stores the two smallest eigenvectors of its
normalized Laplacian.

\begin{figure*}
\subfigure[Spectral graph embedding and the dendrogram for HSC
]{\includegraphics[width=0.43\textwidth]{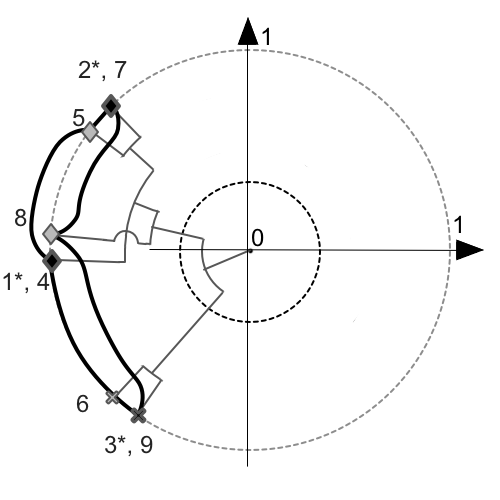}\label{fig:eig9bus}}
\subfigure[Post-islanded network, power flows and islanding performance
metrics]{\includegraphics[width=0.53\textwidth]{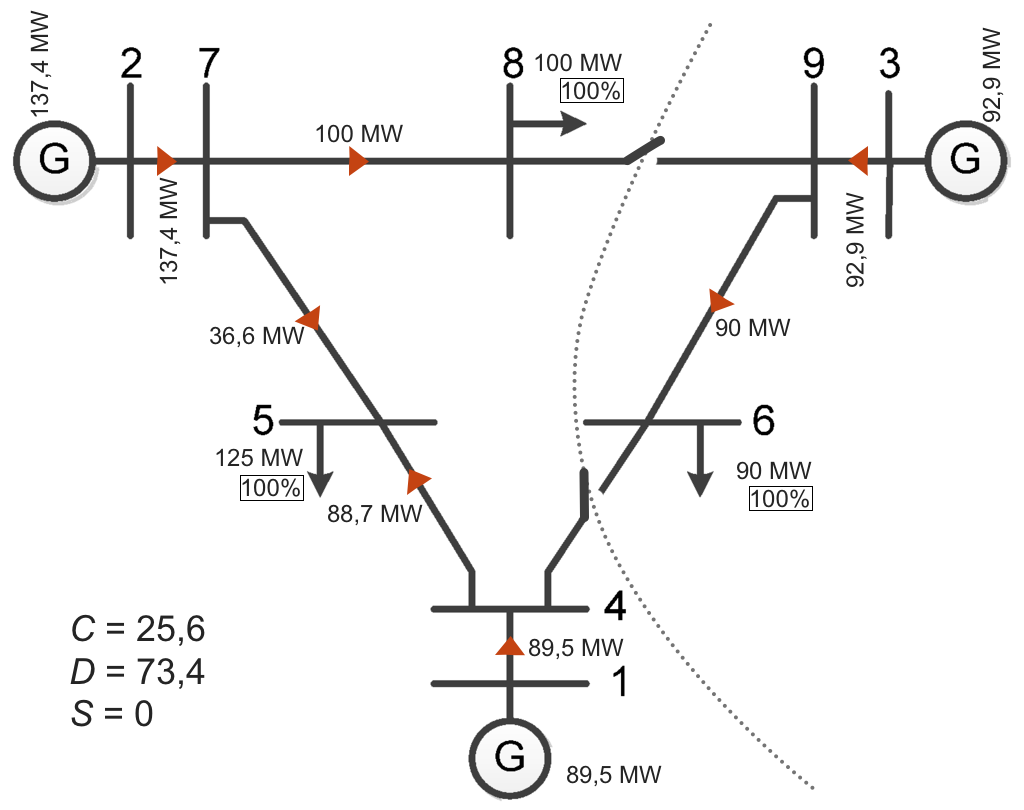}\label{fig:post9bus}}
\caption{An illustration of HSC algorithm}\label{fig:9bus}
\end{figure*}

Each diamond or cross on the unit circle in Figure \ref{fig:eig9bus} represents some
row of matrix $U$ additionally normalized to the unit norm. Each point corresponds
to a graph node (their numbers are listed in the figure with generator nodes marked
with the star), and connecting the points with graph edges (bold curves in the
figure) we obtain the \emph{spectral graph embedding}. Every edge is labelled with a
weight being equal to the distance on the circle between its ends. The dendrogram of
the hierarchical clustering applied to the weighted distance matrix of this graph is
also shown in Figure \ref{fig:eig9bus}. The root of the dendrogram is located at
zero, and the smaller dashed circle in Figure \ref{fig:eig9bus} shows the level, at
which exactly two clusters are separated. The nodes of the first cluster are denoted
with diamonds, while the nodes of the second one are denoted with crosses. It is
worth noting that since we use graph distances, node $1,4^*$ is closer to node $5$
than to node $8$ (which can also be seen at the dendrogram).

The resulting islanding scheme is presented in Figure \ref{fig:post9bus}, which
shows teared lines and equilibrium post-islanding generations, power flows, and
loads (the shares of demand served are shown in frames). Islanding performance
metrics are also presented: generator coherence $C$, disruption $D$ and shed load
$S$ (according to AC-OPF model).

In many applications it is important to keep some pairs of vertices in different
clusters (e.g., generators with low dynamic coupling) and the other pairs should
always be assigned to a single cluster (e.g., highly coupled generators). The
algorithm of constrained spectral clustering \cite{bie2004constrainedclustering}
uses projection matrix technique to generalize the spectral clustering approach to
the case of such must-link and cannot-link constraints, but it can be applied only
for graph bisection (for $K=2$). Constrained spectral clustering is used in
\cite{ding2013CoherencyDisruptionConstrainedSpectralBisection,ding2014constrained,quiros2015constrained}
to control coherent generator groups. In the first step of SCCI algorithm
\cite{ding2013CoherencyDisruptionConstrainedSpectralBisection} the \emph{dynamic
graph} with set $N_g=\{1, ..., n_g\}$ of generators and weights' matrix
$\tilde{\Phi}$ is considered and all generators are divided in two coherent groups
using some spectral bisection algorithm. In the second step unnormalized constrained
spectral clustering is used to select a bisection of grid graph that minimizes flow
disruption and fulfills must-link and cannot-link constraints: generators from the
same coherent group must go to one island while those from different groups cannot
go to one island.

Therefore, the first two eigenvectors $\mathbf{u}_1,\mathbf{u}_2$ of
the generalized eigenvalue problem
$$H^T L(|P|)H \mathbf{u}=\lambda H^T H \mathbf{u}$$
are calculated, where $H$ is a projection matrix. If, without loss
of generality, $s_1' = \{1, ..., n_1\}$, then
$$H = \left(%
\begin{array}{ccc}
  \mathbf{1}_{n_1} & \mathbf{1}_{n_1} & \mathbf{0}_{n_1} \\
  \mathbf{1}_{n_g-n_1} & -\mathbf{1}_{n_g-n_1} & \mathbf{0}_{n_g-n_1} \\
  \mathbf{1}_{n-n_g} & \mathbf{0}_{n-n_g} & I_{(n-n_g)\times (n-n_g)} \\
\end{array}%
\right).$$ After that the rows of matrix
$U=(\mathbf{u}_1,\mathbf{u}_2)$ are split in two clusters using the
$k$-medoids algorithm.
%The resulting bipartition is denoted as
%$\sigma(\mu_2(\tilde{\Phi}|\mathbf{H}),|P|)$, where the first
%argument is a procedure that builds coherent generator groups and
%the second argument is the weights' matrix used to assign loads to
%generator groups.
The same procedure should be applied recursively to the biggest
island of the partition until the desired number of partitions is
obtained.

Another Constrained Spectral Clustering (CSC) methodology of
intentional controlled islanding is proposed in
\cite{quiros2015constrained}. CSC algorithm starts from predefined
groups $s_1^g, ..., s_K^g\subseteq N_g$ of coherent generators
(probably, identified with independent component analysis
\cite{ariff2013coherency} or hierarchical trajectory cluster
analysis \cite{juarez2011characterization}) and aims at forming
islands by distributing loads between these generator groups to
minimize the normalized cut of absolute power flow matrix.

Analogously to HSC algorithm, $K$ first eigenvectors of the
normalized Laplacian of absolute power flow matrix $|P|$ are used to
calculate the graph embedding onto the $K$-dimensional unit sphere.
As in HSC algorithm, the distance between a pair of graph vertices
is evaluated as the length of the shortest path in the embedded
graph (again, cosine distance between incident vertices is
considered). Finally, the $K$-partition of the graph is built by
assigning each vertex to the $i$-th island if and only if the
nearest generator vertex in the embedded graph belongs to group
$s_i^g$. Let us denote with $\sigma_K(|P|~|\mathbf{w},\pi^g)$ a
$K$-partition built by CSC algorithm given absolute power flow
matrix $|P|$, partition $\pi^g=\{s_1^g, ..., s_K^g\}$ of generators
into coherent groups, and vector $\mathbf{w}$ of bus weights.

\begin{figure*}
\subfigure[Spectral graph embedding and load bus assignments for
CSC]{\includegraphics[width=0.43\textwidth]{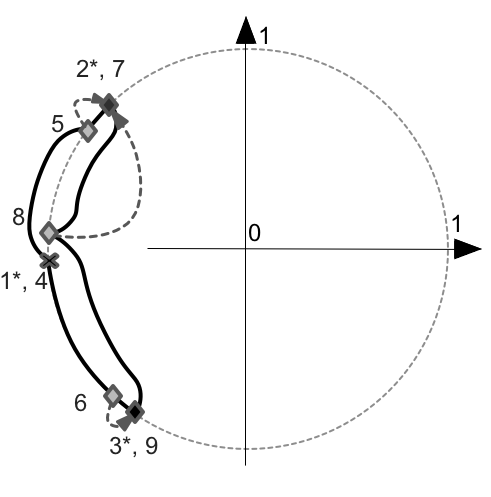}\label{fig:eigcsc9bus}}
\subfigure[Post-islanded network, power flows and islanding performance
metrics]{\includegraphics[width=0.53\textwidth]{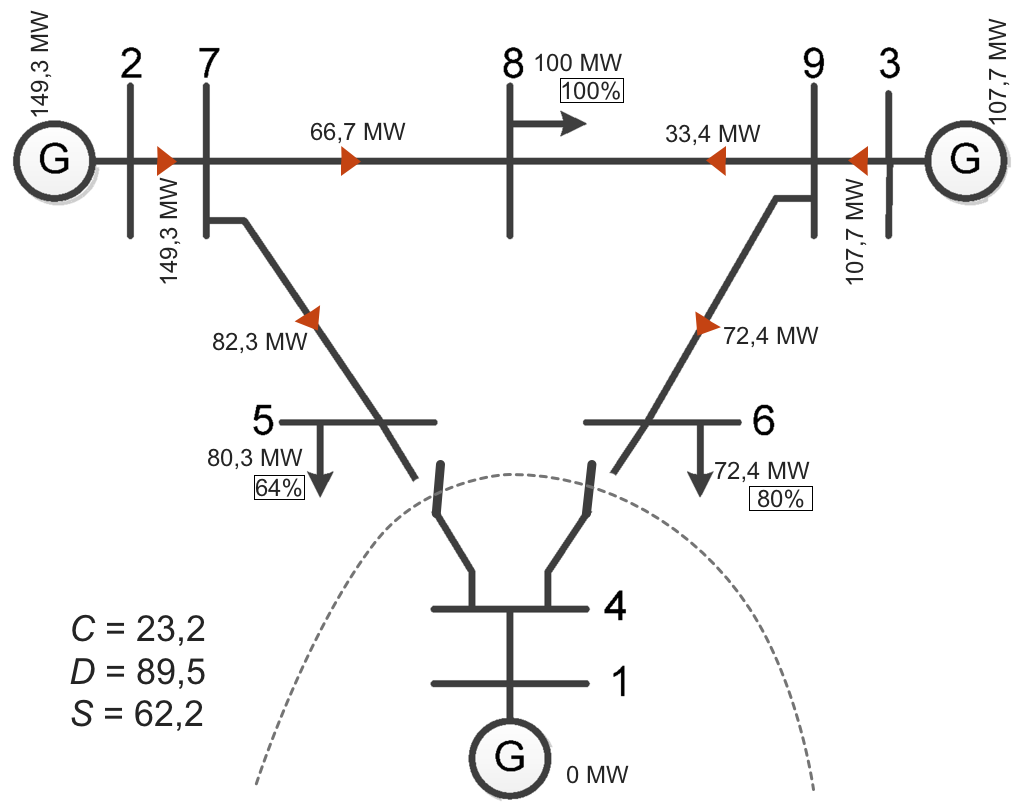}\label{fig:postcsc9bus}}
\caption{An illustration of CSC algorithm}\label{fig:csc9bus}
\end{figure*}

To illustrate CSC algorithm, let us bipartition the 9-bus network (see Figure
\ref{fig:pre9bus}). The spectral graph embedding coincides with that of HSC
algorithm (see Figure \ref{fig:eigcsc9bus}). From Figure \ref{fig:pre9bus} we see
that generators at buses $2$ and $3$ have higher dynamic coupling $\phi_{ij}$, so
let us form the two generator groups: $\{1\}$ (denoted with the black cross in
Figure \ref{fig:eigcsc9bus}) and $\{2,3\}$ (denoted by two black diamonds in Figure
\ref{fig:eigcsc9bus}). Any other bus is assigned to the generator group being
closest to this bus in the graph embedding (see the dashed arrows in Figure
\ref{fig:eigcsc9bus}). Two resulting islands, post-islanding power flows, and
islanding performance metrics are presented in Figure \ref{fig:postcsc9bus}. Load
shedding at buses 5 and 6 is required to balance generation and load in the bigger island.

Like SCCI, CSC algorithm takes care both of generator coherence and power flow
disturbance, but CSC is shown in \cite{quiros2015constrained} to work faster. At the
same time, it disregards shed load (e.g., compare $S$ in Figures \ref{fig:post9bus}
and \ref{fig:postcsc9bus}). Also, under this methodology the size of islands being
created is hardly controllable.

The goal of the present article is to improve spectral clustering algorithms of OCI
(namely, HSC and CSC algorithms) by proposing a methodology, which flexibly accounts
for generator coherence, power flow disruption, shed load, and island sizes.

The deeper insight into spectral clustering techniques can be found in
\cite{von2007tutorial,sanchez2014AdmittanceDisruptionHierSpec,quiros2014determination}
and in the references provided by these articles.

\section{Improved Spectral Clustering Algorithm}\label{sec:Algorithm}

\subsection{Problem Setting}

Effect of islanding criterion on the stability of the islanded grid
was analyzed in
\cite{junchen2015CoherencyDisruptionImbalanceCompare}. An experiment
with IEEE 118-bus scheme has shown that an islanding scheme with
minimum power flow disruption is the most stable one; the scheme
with minimum excess demand results in considerably lower load
shedding at the cost of longer relaxation time, while under a
coherency-based islanding scheme, which ignores disruption and
imbalance, generators fail to stabilize.

It is not clear at the moment to what extent these observations
generalize to other contingency cases and to other power grids:
different performance metrics may be valuable predictors of island
stability in different situations. Therefore, a universal algorithm
of OCI should combine multiple criteria (generator coherency, power
flow disruption, some metric of load shedding, minimum or maximum
island volume, and others) and their relative importance should be
flexibly adjusted.

Distinct to numerous approaches to OCI that consider preserving
generator coherence as a primary optimization goal and calculate
coherent generator groups in advance using the two-time-scale theory
\cite{chow1987Sparsetimescale,chow1988Twotimescale,chow1991Twotimescale}
or recent approaches
\cite{ariff2013coherency,juarez2011characterization}, below the slow
coherency detection technique from
\cite{ding2013CoherencyDisruptionConstrainedSpectralBisection} is
adopted and dynamic coupling $C(\cdot)$ is included into the
optimization criterion along with other performance metrics. OCI is
considered as a multi-objective optimization problem: to find
a $K$-partition $\pi=(s_1,...,s_K)$ of grid graph $\langle
N, E\rangle$ that minimizes the weighted sum of multiple metrics
\begin{equation}\label{eq:F_def}
F(\pi)=\alpha_C C(\pi)+\alpha_D D(\pi)+\alpha_{ECI}
ECI(\pi)+\alpha_S S(\pi)
\end{equation}
and has the limited maximum island volume:
\begin{equation}\label{eq:volume_constraint}
w(s_k)\le W, k=1,...,K.
\end{equation}
The proper choice of weights is discussed in the last section.

Without loss of generality assume $\alpha_S=1$. Then $F(\pi)$ can be
written as
\begin{equation}\label{eq:F_A_S}
F(\pi)=Cut_A(\pi)+S(\pi),
\end{equation}
where
\begin{equation}\label{eq:A_def}
A = (a_{ij})_{i,j=1}^n := \alpha_C \Phi+\alpha_D |P|+\alpha_{ECI} \Delta.
\end{equation}

If shed load is estimated using the maximum flow approximation (see
the previous section), OCI reduces to choosing $x_{ik} \in \{0,1\}$,
$z_{ij} \in \{0,1\}$, $y_{ij}\in \left[-\bar{p}_{ij}z_{ij},
\bar{p}_{ij}z_{ij}\right]$, $l_i \in [0, D_i]$, $g_i\in [0, G_i]$
for $i\in N$, $k=1,...,K$, $ij\in E$ to
\begin{equation}\label{eq:MIQP}
\begin{aligned}
 \text{minimize} ~&\tr X^T L(A)X + \sum_{i=1}^n l_i &\\
\end{aligned}
\end{equation}
subject to constraints:
\begin{align}
\text{island volume} ~& \sum_{i=1}^n w_ix_{ik}\le W
& \forall k=1,...,K, \nonumber\label{eq:volume_constraint2}\\
\text{nodal flow balance} ~& D_i-l_i=g_i+\sum_{j:ji\in
E}y_{ji}-\sum_{j:ij\in E} y_{ij} & \forall i=1,...,n, \nonumber\\
\text{islands' detachment} ~& \sum_{k=1}^K k(x_{ik}-x_{jk}) \le
K(1-z_{ij}), & &\nonumber\\
-&\sum_{k=1}^K k(x_{ik}-x_{jk}) \le K(1-z_{ij}) & \forall ij\in E,\nonumber\\
\text{vertex partitioning} ~& \sum_{k=1}^K x_{ik}=1 & \forall
i=1,...,n.\nonumber
\end{align}

Here $X = (x_{ik})$ is an $n\times K$ indicatory matrix of the grid partition, so
that $x_{ik}=1$ if and only if bus $i$ belongs to island $k$, $y_{ij}$ is a real
power flow through line $ij\in E$ (limited by the maximum flow $\bar{p}_ij$),
$z_{ij}=1$ if and only if line $ij\in E$ lies inside a single island (and, thus, is
not switched off when the grid is partitioned), $l_i$ and $g_i$ are,
correspondingly, shed load volume and real power output at bus $i\in N$, while $D_i$
and $G_i$ being, respectively, real power demand and maximum real power output of
generator at bus $i$.

Since the Laplace matrix is always positively semidefinite, the
problem in hand is convex MIQP with $nK+m$ binary and $2n+m$ real
variables. Nevertheless, numeric optimization packages (e.g., CPLEX)
cannot be applied directly to this problem due to its high
dimension. Below an efficient computational approach is introduced
that avoids the combinatorial explosion.

\subsection{Idea of Algorithm}

Matrix $A$ in expression (\ref{eq:A_def}) is symmetric and
non-negative, so existing algorithms of spectral clustering can
minimize $Cut_A(\pi)$, the first term of cost function
(\ref{eq:F_A_S}). At the same time, distinct to the standard graph
cut problem the sparsity pattern of matrix $A$ does not coincide
with that of graph adjacency matrix due to coupling coefficients
$\phi_{ij}$ that directly ``connect'' non-adjacent generator buses
$i, j\in N_g$. As a result, the classical spectral clustering
algorithm \cite{ng2002spectral} often suggests disconnected islands
with lots of generation-rich and generation-deficient connected
components (and, hence, with poor load shedding). The connectivity
problem can be avoided by using hierarchical spectral clustering
\cite{sanchez2014AdmittanceDisruptionHierSpec} that, in most cases,
generates connected islands.

Absolute power flow matrix $|P|$ is included, among the others, into
matrix $A$, so partition $\pi$ with low $Cut_A(\pi)$ typically has
reasonably low flow disruption $D(\pi)$. In turn, disruption, to
some extent, correlates with load shedding: when losses are
neglected, the inequality
$$S_{EL}(s)=\max\left[Cut_P(s);0\right]\le
Cut_{|P|}(s)=D(s)$$ holds for any island $s\subseteq N$, which means
that low disruption implies low excess load (but not vice versa, in
general). So, when looking for the partition that minimizes a
combination of load shedding and disruption, one can limit attention
to partitions with relatively low disruption.

We use this observation to propose the Improved Spectral Clustering
(ISC) algorithm of controlled islanding. In the first step of the
algorithm a limited set of partitions is obtained with low
$Cut_A(\cdot)$, while in the second step a partition that minimizes
$Cut_A(\cdot)+S_{MF}(\cdot)$ is selected from this set (load
shedding is estimated using the maximum flow model). To detect more
candidate partitions, in the first step different spectral
clustering techniques are combined to build several alternative
graph $rK$-cuts, where $r>0$ is some granularity factor. Then, in
the second step, some of $rK$ islands are merged to obtain a
$K$-partition that minimizes $Cut_A(\cdot)+S_{MF}(\cdot)$.

\begin{figure*}
\subfigure[Initial
plate]{\includegraphics[width=0.3\textwidth]{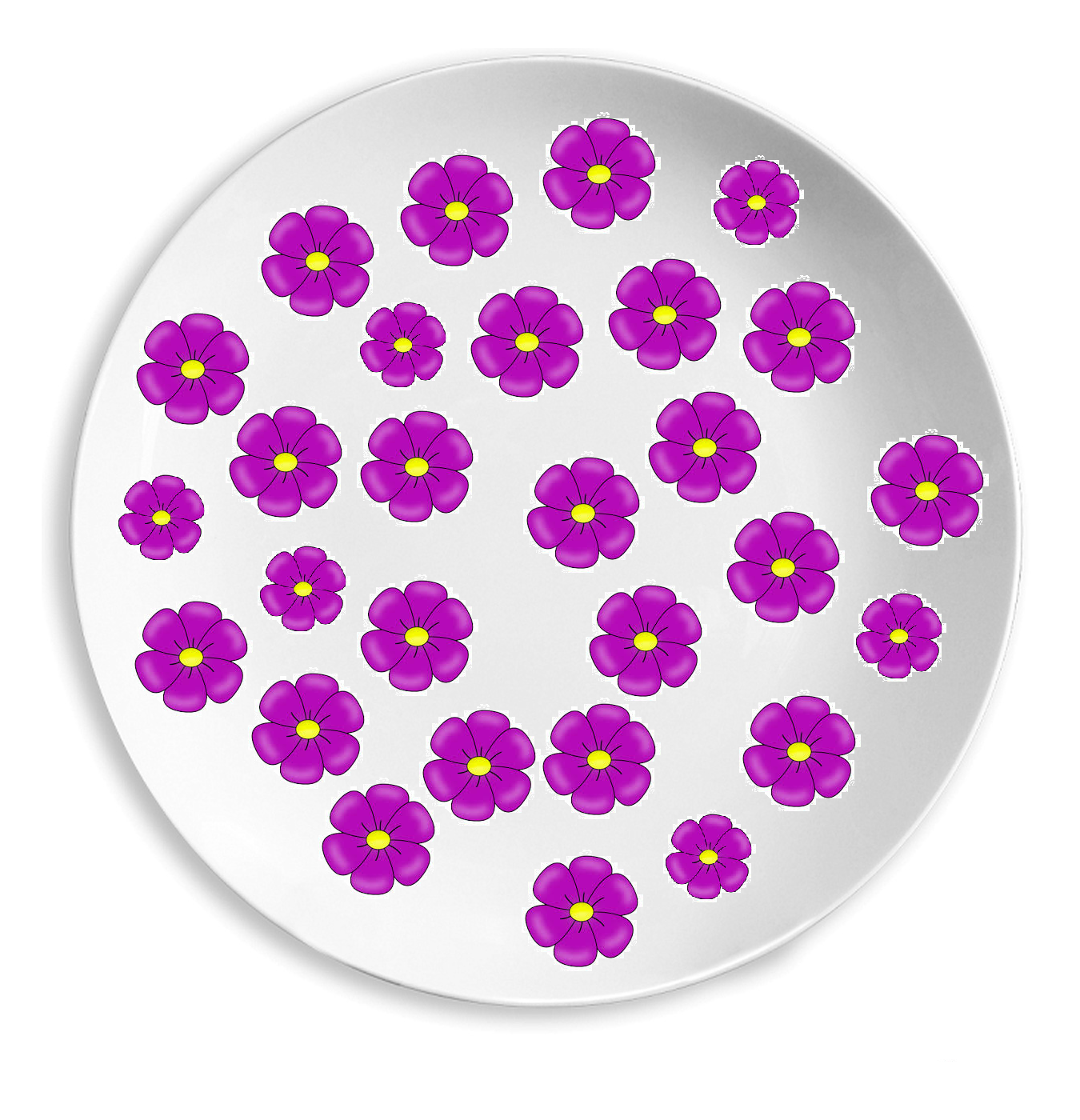}\label{fig:platea}}
\subfigure[Saw
plate]{\includegraphics[width=0.3\textwidth]{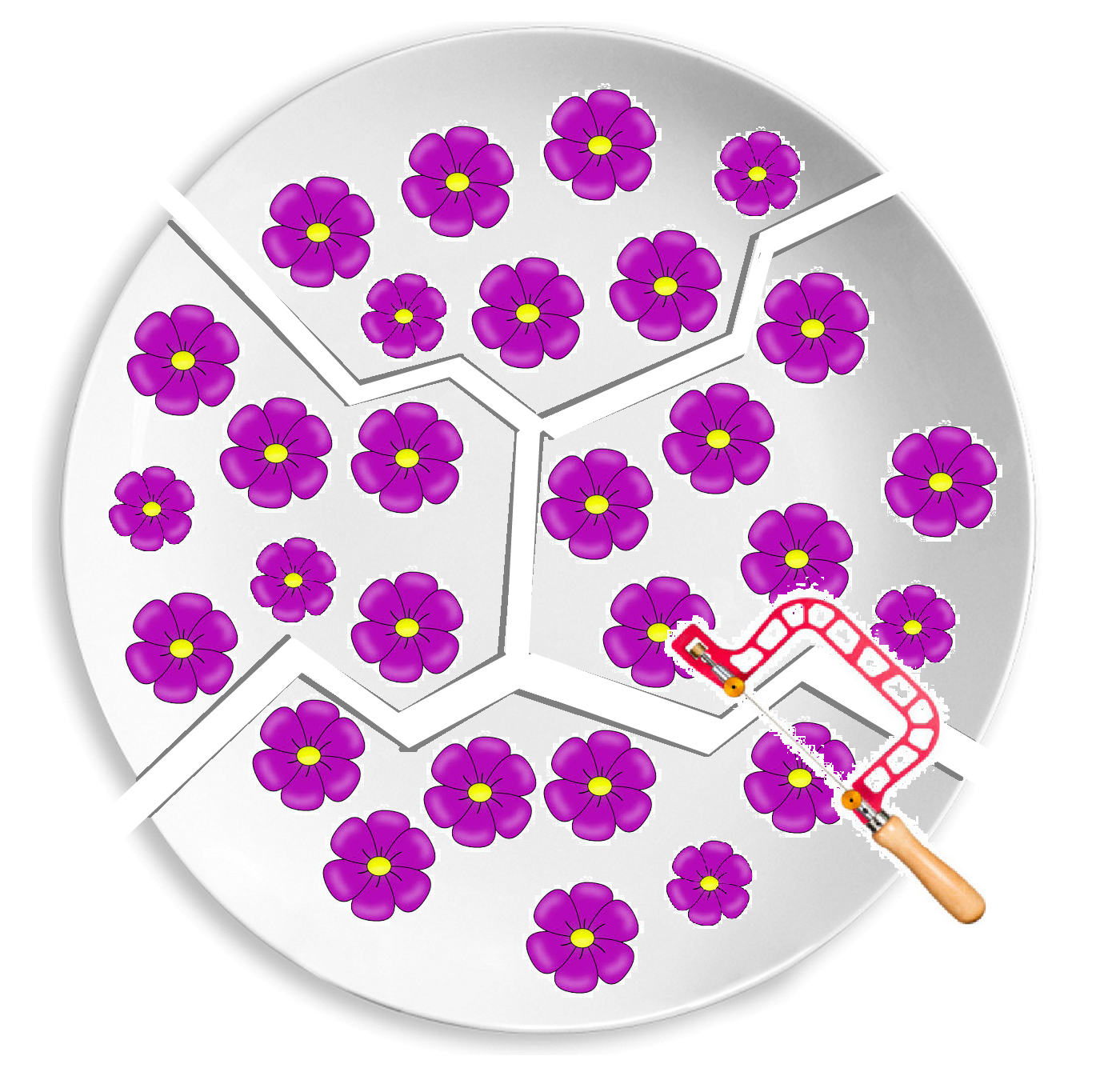}\label{fig:plateb}}
\subfigure[Broken
plate]{\includegraphics[width=0.3\textwidth]{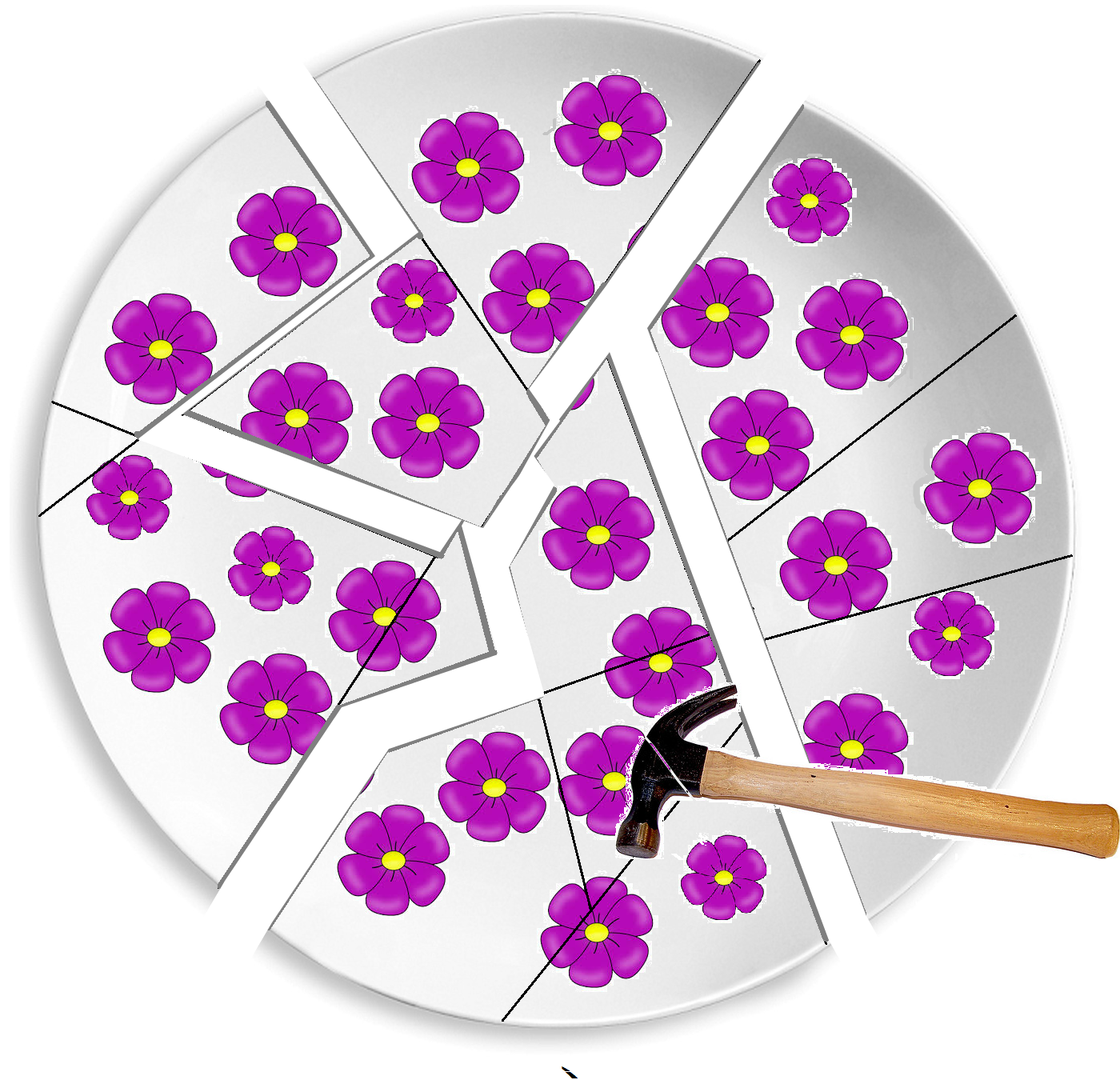}\label{fig:platec}}
\caption{An illustration of the idea of the
algorithm}\label{fig:plate}
\end{figure*}

The following story illustrates the idea of the algorithm. Imagine
someone has a porcelain plate with flowers painted on it (see
Fig.~\ref{fig:platea} and wants to divide it into four pieces of
roughly equal size keeping flowers unbroken. To obtain such pieces a
plate can be sawed carefully (one of possible solutions is shown in
Fig.~\ref{fig:plateb} but it is extremely time-consuming. Instead,
one can break a plate into small pieces with a strong hammer blow
and then glue some pieces back to compose as much entire flowers as
possible. Although the latter approach may result in a suboptimal
solution (e.g., five flowers are broken in Fig.~\ref{fig:platec}),
it is much faster and can be the only alternative when time is
expensive.

Details of the algorithm are explained in the next two subsections.

\subsection{Step 1: Spectral Clustering}\label{seq:spectral}

Below we simplify presentation by assuming that $\alpha_{ECI}=0$ in
(\ref{eq:F_def}). To apply normalized spectral clustering, the
balanced cut problem with maximum cluster volume constraints is
replaced with the corresponding minimum normalized cut problem
without cluster volume constraints. Seven different strategies are
used in parallel in the first step to build overdetailed partitions.
Then, in the second step of the algorithm, each of these partitions
is rolled up into a $K$-partition.

\begin{enumerate}[I]
    \item\label{strategy_fixedA}\textbf{Fixed-granularity strategy for matrix $A$.}
    Granularity factor $r_1>1$ (which is a tunable
parameter of the algorithm) is chosen and $r_1K$-partition
$\pi^\text{\ref{strategy_fixedA}} := \chi_{r_1K}(A|\mathbf{w})$
    of the grid is built with HSC algorithm \cite{sanchez2014AdmittanceDisruptionHierSpec}.
%Cluster volume constraints (\ref{eq:volume_constraint}) are seldom
%violated for $r_1 \ge 2$ and the practical range of $W$: when it
%occasionally happens, $r_1$ is increased and a more granulated
%partition is calculated.
    \item\label{strategy_fixedP}\textbf{Fixed-granularity strategy for matrix $|P|$.}
    It has already been noted that low disruption implies low excess load,
    so partitions with low $D(\cdot)$ seem more likely to minimize
    cost function (\ref{eq:F_def}) than those with low $C(\cdot)$.
    Therefore, partition $\pi^\text{\ref{strategy_fixedP}}:=\chi_{r_2K}(|P|~|\mathbf{w})$ is computed, where granularity factor
    $r_2>1$ is a yet another tunable parameter.
    \item \label{strategy_minA}\textbf{Minimum-granularity strategy for matrix $A$.}
    From inequality (\ref{ineq:spectral}) it follows that exactly
    $K$ Laplacian eigenvectors are
    enough to characterize a $K$-partition with small normalized
    cut. So, overdetailed partitions
    in Strategies \ref{strategy_fixedA} and \ref{strategy_fixedP}
    are forced by the need to satisfy volume
    constraints (\ref{eq:volume_constraint}) and to leave some combinatorial space for the shed
    load optimization. Instead, in Strategy \ref{strategy_minA} a
    partition $\pi^\text{\ref{strategy_minA}} := \chi_{K'}(A|\mathbf{w})$ is chosen with minimum granularity $K'\in \{K, ..., r_1K\}$ that
    satisfies cluster volume constraints (\ref{eq:volume_constraint}) (i.e., if $s\in \pi^\text{\ref{strategy_minA}}$ then $w(s)\le W$).
    \item \label{strategy_minP}\textbf{Minimum-granularity strategy for matrix $|P|$.}
    Similarly to the previous strategy, we choose the partition $\pi^\text{\ref{strategy_minP}} := \chi_{K''}(|P|~|\mathbf{w})$ with minimum granularity $K''\in \{K, ..., r_2K\}$ that
    satisfies cluster volume constraints (\ref{eq:volume_constraint}).
    \item \label{strategy_CSC}\textbf{Minimum-granularity-refined CSC algorithm.}
    CSC algorithm minimizes disruption under given groups of coherent generators. In
our methodology selection of coherent generator groups on the basis
of dynamic coupling $\Phi$ (or alternative generator coherence
metric) is a part of the optimization process.
    In Strategy \ref{strategy_CSC} coherent generator
groups for CSC algorithm are calculated by HSC algorithm applied to
reduced dynamic coupling matrix $\tilde{\Phi}$. Consequently,
$\pi^\text{\ref{strategy_CSC}} :=\sigma_{K'''}(|P|~|
\chi_{K'''}(\tilde{\Phi}|\mathbf{w}), \mathbf{w})$ for the minimal
granularity $K'''\in \{K, r_3K\}$ that allows to satisfy island
volume constraints (\ref{eq:volume_constraint}) (as before, $r_3>0$
is a tunable parameter).
    \item \label{strategy_fixedPA}\textbf{Fixed-granularity sequential strategy for $|P|$ and $A$.}
    HSC algorithm generates a partition $\chi_K(|P|~|\mathbf{w})$ that has low disruption but neglects generator coherency.
    So, recursive bisection for graph edge weights' matrix $A$ is applied to partition
$\chi_K(|P|~|\mathbf{w})$ to obtain more granulated $r_4K$-partition
(again, $r_4>1$ is a tunable parameter),
    which combines low disruption and coherency.

Let us denote with $\nu(\sigma, \pi|\mathbf{w})$ a partition
obtained from partition $\pi$ by splitting island $s\in \pi$, which
has the biggest volume $w(s)$, with a bisection procedure $\sigma$.
Define
$$\nu_K(\sigma, \pi|\mathbf{w}) = \nu(\sigma, \nu(...,\nu(\sigma, \pi|\mathbf{w})| ...)|\mathbf{w})$$
a partition being a result of $K$ recursive bisections $\nu$ of the
initial partition $\pi$.
    Then, $\pi^\text{\ref{strategy_fixedPA}} :=
\nu_{(r_4-1)K}(\chi_2(A|\mathbf{w}),
    \chi_K(|P|~|\mathbf{w})|\mathbf{w})$. The idea behind this strategy
    is that we never miss low-disruption
    partition $\chi_K(|P|~|\mathbf{w})$ and potentially can improve
    by joining islands in another order in the second step of the algorithm.
    \item \label{strategy_meet}\textbf{Crossing CSC and HSC partitions.}
    CSC algorithm builds a partition $\sigma_K(|P|~|
\chi_K(\tilde{\Phi}|\mathbf{w}), \mathbf{w})$ focusing mainly on
generator coherence. On the contrary, HSC algorithm constructs a
    partition $\chi_K(|P|~|\mathbf{w})$ caring only for flow disruption.
    A finer partition $\pi^\text{\ref{strategy_meet}} := \sigma_K(|P|~|
\chi_K(\tilde{\Phi}|\mathbf{w}), \mathbf{w}) \wedge
\chi_K(|P|~|\mathbf{w})$
    is obtained as a \emph{meet} of these two partitions in the \emph{aggregation
    lattice}. Set $s\subseteq N$ belongs to the meet
    if an only if $s = s_1 \cap s_2$ for some $s_1\in
    \sigma_K(|P|~|
\chi_K(\tilde{\Phi}|\mathbf{w}), \mathbf{w})$ and $s_2 \in
    \chi_K(|P|~|\mathbf{w})$.
\end{enumerate}
We limit ourselves to the above seven strategies, although other
approaches to granulated partition construction are also possible.

\subsection{Step 2: Partitioning Aggregated Grid}\label{sec:MIQP}

In the second step of the algorithm the granulated partition is
transformed into $K$-partition by fusing some islands together. Each
of detailed partitions $\pi^\text{\ref{strategy_fixedA}}, ...
\pi^\text{\ref{strategy_meet}}$ is processed separately and,
probably, in parallel.

First of all, all connected components of islands in partition
$\pi^i$ are detached making separate islands. Let $\pi=(s_1, ...,
s_{n'})$ be the resulting detailed partition where each island is
connected ($n'>K$).

An aggregated grid is built such that each island $s_k$ in $\pi$
becomes its vertex with:
$$
\begin{array}{rl}
    \text{maximum real power output} & G_k':=\sum_{i\in s_k}G_i, \\
    \text{current real power output} & g_k':=\sum_{i\in s_k}g_i, \\
    \text{real power demand} & D_k':=\sum_{i\in s_k}D_i, \\
    \text{current real power load} & d_k':=\sum_{i\in s_k}d_i, \\
    \text{injection} & p_k':=p(s_k), \\
    \text{volume} & w_k':=w(s_k). \\
\end{array}
$$

Also, for all island pairs $k,k'=1, ..., n'$ define
$$
\begin{array}{rl}
    \text{aggregated real power flows} & p_{kk'}'=\sum_{i\in s_k,j\in
s_{k'}}p_{ij}, \\
    \text{real power limits} & \bar{p}_{kk'}'=\sum_{i\in s_k,j\in
s_{k'}}\bar{p}_{ij}, \\
    \text{dynamic coupling coefficients} & \phi_{kk'}'=\sum_{i\in s_k,j\in s_{k'}}\phi_{ij}. \\
\end{array}%
$$

Edge $kk'$ is included into edge set $E'$ of the aggregated grid if
$\bar{p}_{kk'}'>0$. Let $m':= |E'|$ be the edge count of the
aggregated grid, $m'\le n'(n'-1)/2$.

OCI problem (\ref{eq:MIQP}) for the aggregated grid is a relaxation of the same OCI
problem for the original power grid, because it replaces a group of detailed
constraints for vertices of one island with a lump-sum constraint for the island as
a whole. The aggregated grid problem is MIQP with $Kn'+m'$ binary and $2n'+m'$ real
variables. Its dimension depends on the dimension of the aggregated grid but not on
the dimension of the original grid, and if granularity $n'$ of partition $\pi$ is
small enough, this problem can be solved in eligible time using exact algorithms
implemented in commercial numeric solvers (we use CPLEX 12.6.2.0).\footnote{For
further acceleration of calculations the problem in hand can be reduced to the
mixed-integer linear problem (MILP) with the techniques described in
\cite{fanpardalos2012SLDC_MILP,trodden2013milp} but in the present article this
possibility is not studied in detail.}

If there are too many disconnected islands in $\pi^i$, dimension
$n'$ of the aggregated grid can still be too high for exact
algorithms. In this case we suggest decreasing its dimension with
the following greedy heuristics.

At each iteration this heuristics simplifies the partition by
joining a pair of adjacent islands to minimize cost function
(\ref{eq:F_A_S}) while fulfilling island volume constraints. The
shed load in (\ref{eq:F_A_S}) is estimated with the excess demand
$S_{ED}(\cdot)$:
\begin{equation}\label{eq:F_simpl}
F_{ED}(\pi):=Cut_A(\pi)+S_{ED}(\pi)=Cut_A(\pi)+\sum_{s\in \pi}
\max[p(s),0].
\end{equation}

The pseudo code is presented in Listing \ref{lst:greedy}.
\begin{algorithm}[!ht]
\floatname{algorithm}{Listing} \caption{Greedy heuristic algorithm
for OCI}\label{lst:greedy}
\begin{algorithmic}[1]
\algblockdefx{With}{EndWith} [1]{\textbf{with} #1}{\textbf{end}}
\Function{GreedyPartition}{$Grid$, $K$} \Comment{cuts graph $Grid$
into $K$ parts} \With{Grid}\State $\pi\gets \bigcup_{i\in N}\{i\}$
\Comment{start from the finest partition} \EndWith
\While {$|\pi|>K$}
\Comment{until requested island count is reached} \State
    $\pi \gets \Argmin\{ F_{ED}(\pi')~|\pi' \in \Call{Coarsen}{\pi}, \pi' = \Call{ConnectedComponents}{\pi'}\}$
\EndWhile \State \Return{$\pi$} \EndFunction

\State

\Function{Coarsen}{$\pi$} \Comment{all admissible partitions
obtained by fusing two islands in $\pi$} \State
\Return{$\{\pi'~|~\pi'=\{\{s\cup
s'\}\}\cup(\pi\backslash\{s\}\backslash \{s'\}) \text{ for some
 }s,s'\in \pi,w(s\cup
s')\le W\}$ }
 \EndFunction
\end{algorithmic}
\end{algorithm}

The $K$-partition calculated by this greedy heuristics is also used
as a record in the exact partitioning algorithm of the aggregated
grid.

Finally, ISC algorithm selects the partition with the lowest cost of
MIQP solution for seven aggregated grid built on the basis of
partitions $\pi^\text{\ref{strategy_fixedA}}, ...,
\pi^\text{\ref{strategy_meet}}$.

The pseudo code of ISC algorithm is presented in
Listing~\ref{lst:twostep}.

\begin{algorithm}[!ht]
\floatname{algorithm}{Listing} \caption{Improved Spectral Clustering
algorithm for OCI}\label{lst:twostep}
\begin{algorithmic}[1]
\Function{ImprovedSpectralPartition}{$Graph, K$}
\algblockdefx{With}{EndWith} [1]{\textbf{with} #1}{\textbf{end}}

\With{$Graph$} \Comment{assignments can be done in parallel} \State
$\pi^\text{\ref{strategy_fixedA}} \gets \chi_{r_1K}(A|\mathbf{w})$
\State $\pi^\text{\ref{strategy_fixedP}} \gets
\chi_{r_2K}(|P|~|\mathbf{w})$ \State $K'\gets \min\left\{k=
K,...,r_1K: s\in \chi_k(A|\mathbf{w}) \Rightarrow w(s)\le W\right\}$
\State $\pi^\text{\ref{strategy_minA}} \gets
\chi_{K'}(A|\mathbf{w})$ \State $K''\gets \min\left\{k= K,...,r_2K:
s\in \chi_k(|P|~|\mathbf{w}) \Rightarrow w(s)\le W\right\}$ \State
$\pi^\text{\ref{strategy_minP}} \gets \chi_{K''}(|P|~|\mathbf{w})$
\State $K'''\gets \min\left\{k= K,...,r_3K: s\in \sigma_k(|P|~|
\chi_k(\tilde{\Phi}|\mathbf{w}), \mathbf{w}) \Rightarrow w(s)\le
W\right\}$ \State $\pi^\text{\ref{strategy_CSC}} \gets
\sigma_{K'''}(|P|~| \chi_{K'''}(\tilde{\Phi}|\mathbf{w}),
\mathbf{w})$ \State $\pi^\text{\ref{strategy_fixedPA}} \gets
\nu_{(r_4-1)K}(\chi_2(A|\mathbf{w}),
    \chi_K(|P|~|\mathbf{w})|\mathbf{w})$
\State $\pi^\text{\ref{strategy_meet}} \gets \sigma_K(|P|~|
\chi_K(\tilde{\Phi}|\mathbf{w}), \mathbf{w}) \wedge
\chi_K(|P|~|\mathbf{w})$
    \EndWith \For
{${i=\text{\ref{strategy_fixedA}}, ...,
\text{\ref{strategy_meet}}}$} \Comment{iterations can be performed
in parallel}
    \State $\pi \gets \Call{ConnectedComponents}{\pi^i}$
    \State $Graph' \gets \Call{AggregatedGraph}{Graph,\pi}$
    \If {$|\pi|\ge K^{max}$} \Comment{if partition dimension is too high}
        \State $\pi \gets \Call{GreedyPartition}{Graph',
        K^{max}}$ \Comment{decrease dimension}
      \State $Graph' \gets \Call{AggregatedGraph}{Graph',\pi}$
 \EndIf
    \State $\pi_0 \gets \Call{GreedyPartition}{Graph', K}$ \Comment{record for exact algorithm}
    \State $[\pi', Cost^i] \gets \Call{SolveMIQP}{Graph', K, \pi_0}$ \Comment{exact solution with MIQP solver}
    \State $\pi^{i*} \gets \Call{FullPartition}{\pi'}$ \Comment{obtain $K$-partition of $Grid$ graph}
 \EndFor

\State
$j:=\Argmin_{i=\text{\ref{strategy_fixedA}},...,\text{\ref{strategy_meet}}}
Cost^i$ \Comment{the best calculated $K$-partition} \State \Return
{$\pi^{j*}$} \EndFunction
\end{algorithmic}
\end{algorithm}

\subsection{Numeric Example}\label{sec:Example}

\begin{figure*}
\subfigure[Strategy I (for matrix
$A$)]{\includegraphics[width=0.47\textwidth]{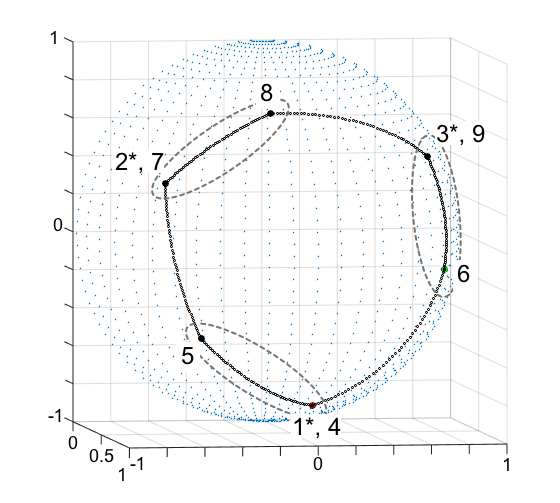}\label{fig:eig2A9bus}}
\subfigure[Strategy II (for matrix
$|P|$)]{\includegraphics[width=0.47\textwidth]{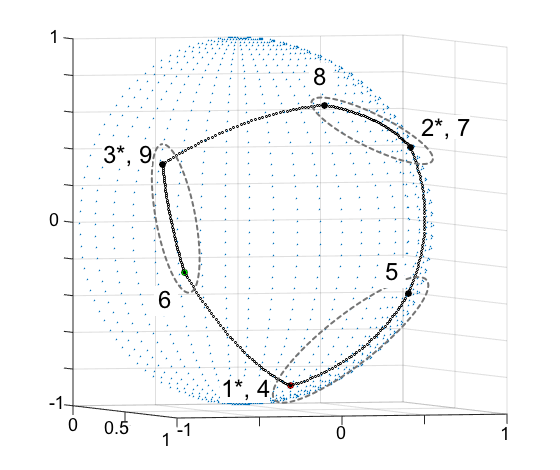}\label{fig:eig2P9bus}}
\caption{Spectral graph embedding and 3-cluster partitioning for Strategies I and
II}\label{fig:9bus2}
\end{figure*}

This subsection illustrates the proposed ISC algorithm by partitioning the 9-bus
network (see Figure \ref{fig:pre9bus}) to minimize the sum of dynamic coupling, flow
disruption, and shed load (with unit weights). The desired number of clusters $K=2$
and granularity factor is set to $r_i=1.5, i=1,...,4$, so, strategies I, II, and VI
in the first step of ISC algorithm try to partition the network into $3$ islands.
Therefore, three eigenvectors are calculated, and the spectral graph embedding in
Strategies I and II is three-dimensional. Strategy I is based on matrix $A$ (a
linear combination of dynamic coupling matrix $\Phi$ and absolute power flow matrix
$|P|$), while Strategy II relies solely on matrix $|P|$. Spectral embeddings for
these strategies (see Figures \ref{fig:eig2A9bus} and \ref{fig:eig2P9bus}
respectively) differ a bit, but hierarchical clustering of the nodes of the spectral
embedding results in the same three islands $\{1,4,5\}, \{2,7,8\}$, and $\{3,6,9\}$
outlined in Figures \ref{fig:eig2A9bus} and \ref{fig:eig2P9bus}.

In the second step of the algorithm, MIQP (\ref{eq:MIQP}) is solved for the 3-nodal
aggregated network, and the resulting bipartition is equal to that calculated by HSC
algorithm (see Figure \ref{fig:post9bus}). In the considered example, Strategy IV
(based on HSC algorithm) and Strategy V (based on CSC algorithm) work as explained
in Section \ref{sec:Spectral} above, and the resulting islanding schemes are shown
in Figures \ref{fig:post9bus} and \ref{fig:postcsc9bus} respectively.

\begin{figure*}
\subfigure[Spectral graph embedding and the
dendrogram]{\includegraphics[width=0.47\textwidth]{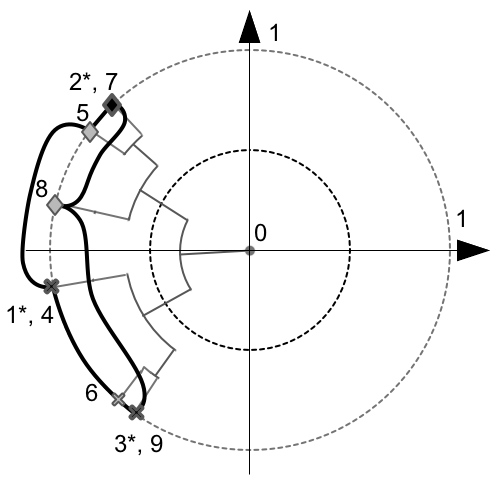}\label{fig:eigA9bus}}
\subfigure[Post-islanded network, power flows and islanding performance
metrics]{\includegraphics[width=0.47\textwidth]{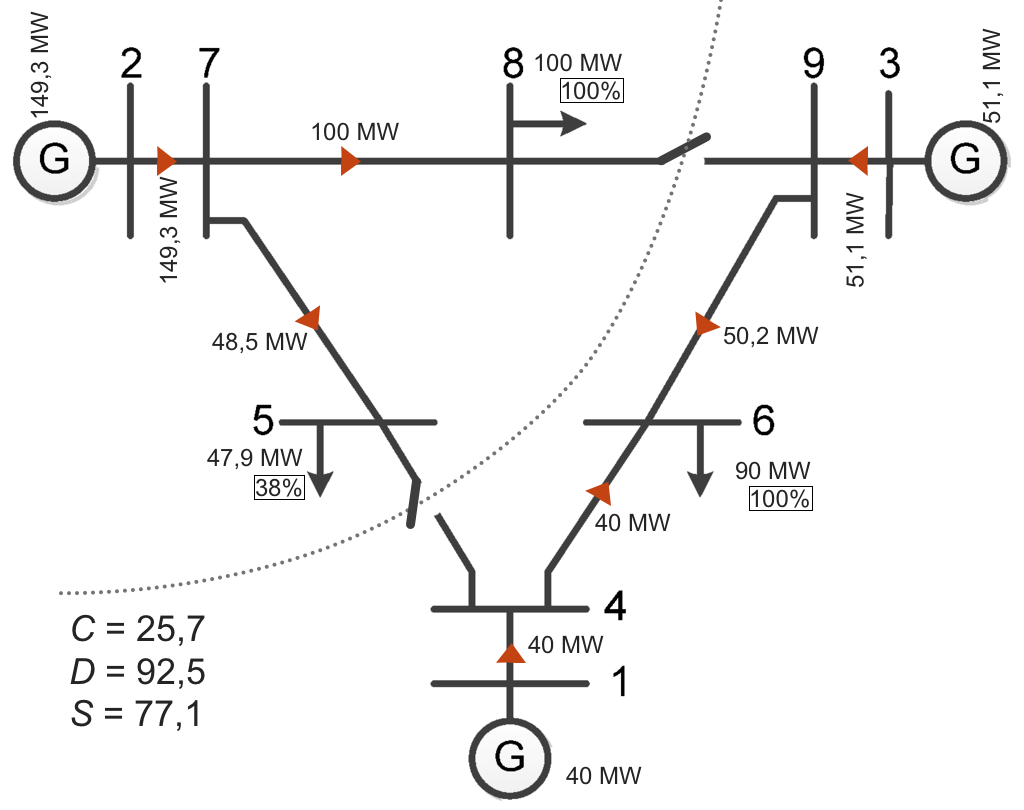}\label{fig:postA9bus}}
\caption{Steps of ISC algorithm: Strategy~III}\label{fig:9bus3}
\end{figure*}

The two-dimensional spectral embedding for Strategy III based on matrix $A$ is
presented in Figure \ref{fig:eigA9bus}. The corresponding islanding scheme,
post-islanding power flows, and islanding performance metrics are shown in Figure
\ref{fig:postA9bus}.

\begin{figure*}
\subfigure[Strategy VI: using CSC algorithm for sequential granulation of
$\chi_K(|P|~|\mathbf{w})$
]{\includegraphics[width=0.47\textwidth]{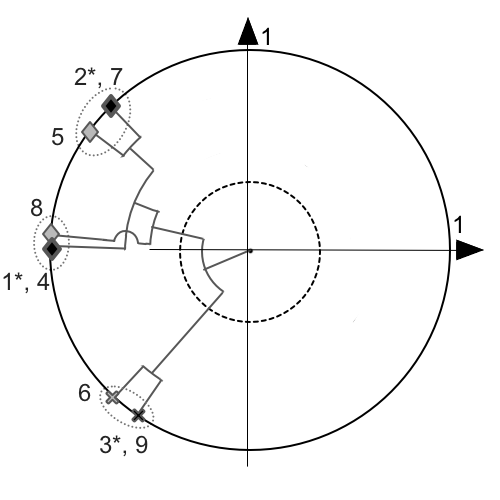}\label{fig:eigseq9bus}}
\subfigure[Strategy VII: the meet of partitions generated by HSC and CSC
]{\includegraphics[width=0.47\textwidth]{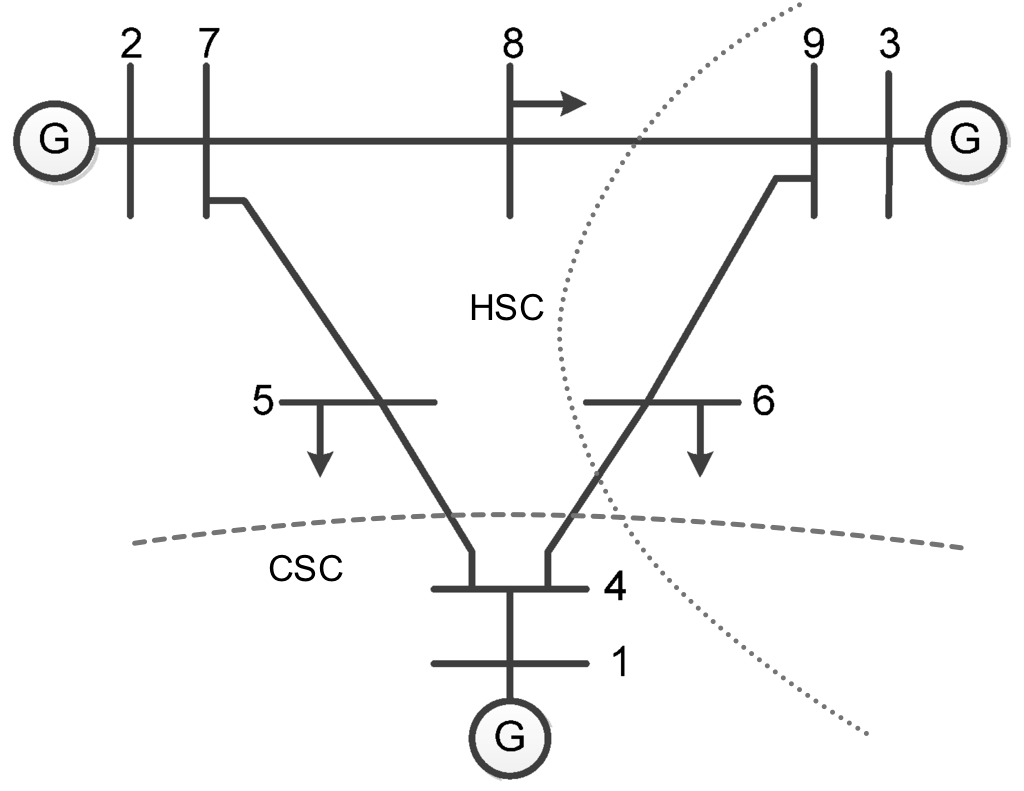}\label{fig:cross9bus}} \caption{ISC
algorithm: Strategies VI and VII}\label{fig:9bus4}
\end{figure*}

Strategy VI starts with the $K$-partition of the network generated by HSC algorithm
(see Figure \ref{fig:eig9bus}) and sequentially bisects the biggest island using CSC
algorithm until $r_4K$ islands are obtained. The three resulting islands are
encircled by dotted lines in Figure \ref{fig:eigseq9bus}. In Strategy VII we meet
partitions generated by HSC algorithm (the dotted line in Figure
\ref{fig:cross9bus}) and CSC algorithm (the dashed line in Figure
\ref{fig:cross9bus}) and obtain three islands: $\{1,4\}$, $\{2,5,3,8\}$, and
$\{3,6,9\}$. In the second step of the algorithm, MIQP (\ref{eq:MIQP}) is solved for
the corresponding aggregated 3-nodal networks, and the resulting partitions coincide
with that shown in Figure \ref{fig:post9bus} both for Strategies VI and VII.

Finally, the islanding scheme with the minimum cost is selected from the three
distinct schemes (those shown in Figures \ref{fig:post9bus}, \ref{fig:postcsc9bus},
and \ref{fig:postA9bus}) obtained by running two steps of ICS algorithm for
Strategies I-VII.

\section{Performance Evaluation}\label{sec:Experiments}

\subsection{Experimental Setup}

Three test systems from MATPOWER 5.1 Simulation Package library
\cite{zimmerman2011matpower} were used to evaluate performance of
the proposed ISC algorithm (see Table \ref{tab:testgrids}). These
relatively big systems was taken to verify computational efficiency
of ISC algorithm and its applicability to bulky real-world grids.
\begin{table}[!ht]
\caption{Test power systems}
\label{tab:testgrids}
\begin{tabular}{llccc}
\hline\noalign{\smallskip}
  Notation & Power system & Buses & Generators & Lines \\
\noalign{\smallskip}\hline\noalign{\smallskip}
  SMALL & IEEE 118-bus test system \cite{IEEE118} & 118 & 54 & 186 \\
  MEDIUM & Polish system, winter peak 1999-2000 \cite{case2383wp} & 2383 & 327 & 2896 \\
  LARGE & European system (PEGASE project) \cite{fliscounakis2013contingency} & 9241 & 1445 & 16049 \\
\noalign{\smallskip}\hline
\end{tabular}
\end{table}

For each of three considered power systems 100 test cases were
generated by switching off 5 random generators, breaking 5 random
lines, and multiplying all demands by a random factor from 1 to 2.
For every case consistent power flows were calculated by solving
AC-OLS problem with \texttt{opf} routine of MATPOWER 5.1
\cite{zimmerman2011matpower} for dispatchable loads.\footnote{To
some extent these starting conditions can be interpreted as a
situation in contingency case after the load shedding program run to
balance demands and available generation/transmission capacities. It
is assumed that these efforts where not enough to stabilize the
system, and controlled islanding is performed.}

\begin{figure*}
\subfigure[Density of power flow from bus 37 to
34]{\includegraphics[width=0.5\textwidth]{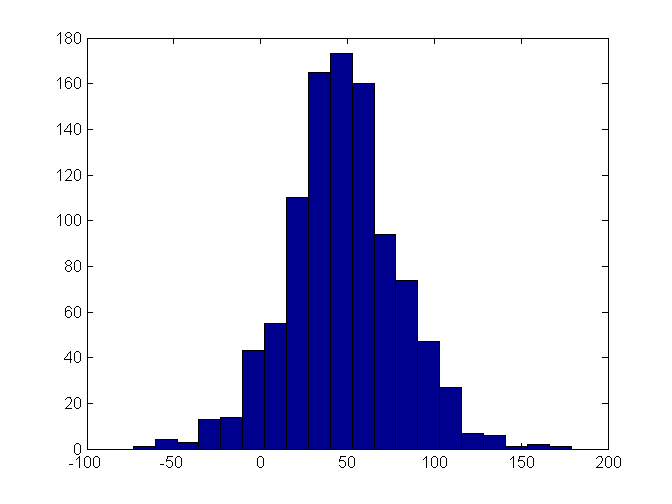}\label{fig:flow}}
\subfigure[Density of real power output at bus
12]{\includegraphics[width=0.5\textwidth]{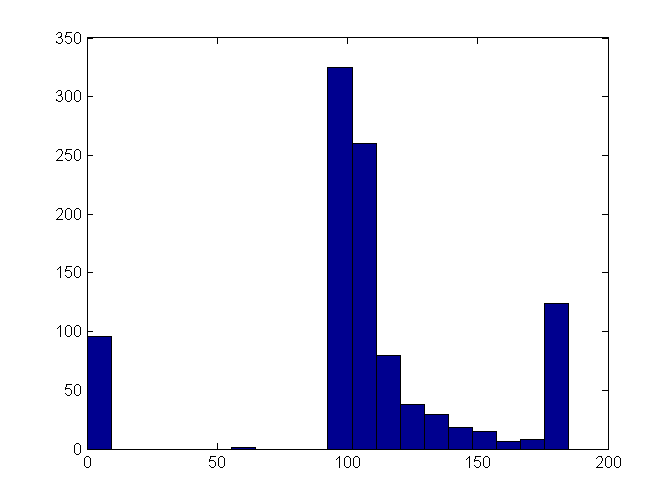}\label{fig:gen}}
\caption{Examples of power flow variables' variation for 100 cases
in SMALL system}\label{fig:density}
\end{figure*}

These power flows and corresponding generators' angles were used to
calculate stability indicators (power flow disruption $D(\cdot)$ and
dynamic coupling $C(\cdot)$) when an islanding operation is planned.
Information on generator inertia constants was not available, so
equal inertia constants were assumed when coefficients of dynamic
coupling $\phi_{ij}$ were calculated.

The sample line power flow and the sample real power output in SMALL
grid are depicted in Figure \ref{fig:density}. High variance of
these variables shows that considered collection of cases represents
a wide range of power flow conditions in a grid. Such massive
durability testing of OCI algorithms under the broad variety of grid
and flow conditions is rarely performed in the existing literature.
The only known exception is \cite{hamon2015two}, where partitioning
algorithms were tested under 2000 different power flow conditions in
SMALL system.

We did not perform transient stability analysis to verify system
stability after islanding, as the main concern of ISC algorithm was
to improve the partition quality for existing island stability
indicators and their combinations.\footnote{Some analysis of
different islanding performance metrics can be found in
\cite{quiros2015constrained}.} At the same time, although simplified
metrics of load shedding (such as excess demand) were employed when
an optimal islanding scheme was searched, for performance evaluation
the AC-OLS model was applied to the system partitioned according to
seven partitioning strategies of ISC algorithm.
%Such an approach
%allows to obtain the better estimate of load shedding in case of
%CSC, HSC, and ISC algorithms correspondingly.

Additional constraints were imposed when AC-OLS problem was solved
for the islanded system: generators' output was limited by the
short-term ramp rate and previously shed load could not be restored
in the process of islanding. So, controlled islanding always results
in extra load shedding compared to the pre-islanding system state.

Dynamic coupling, power flow disruption, and shed load had equal
weights $\alpha_C=\alpha_D=\alpha_S=1$ in the optimization criterion
(\ref{eq:F_def}) of ISC algorithm. The requested island count $K=4$
was selected for all three systems and granularity factors were set
to $r_1=...=r_4=4$, i.e., 16 islands were demanded when calculating
detailed partitions. The maximum island volume was
$W=\frac{3}{8}W(N)$. In particular, this means that any admissible
partition with four connected islands has at least three big islands
(those having the volume $\ge \frac{W(N)}{4}$).

\subsection{Partitioning Quality}

ISC is an approximate algorithm, and there are several possible
sources of its inaccuracy that should be inspected.

Firstly, there may be the discrepancy between real objectives of controlled
islanding (i.e., preserving the transient stability and preventing a blackout with
minimum load shedding) and their representation in the model (\ref{eq:F_def}) (a
combination of computationally efficient metrics). Although being critical for the
final efficiency of an islanding technique, analysis of model adequacy falls beyond
the scope of the present article, which concentrates on the existing OCI metrics.

Another important aspect is the method used to evaluate load
shedding. The estimate $S_{AC}(\cdot)$ of the shed load volume
calculated from AC-OLS model is assumed the most accurate one, but
ISC algorithm employs its maximum-flow relaxation $S_{MF}(\cdot)$ to
reduce the problem to MIQP (\ref{eq:MIQP}). Large discrepancy
between these metrics may sufficiently distort the optimization
criterion.

\begin{figure*}
\subfigure[SMALL
system]{\includegraphics[width=0.49\linewidth]{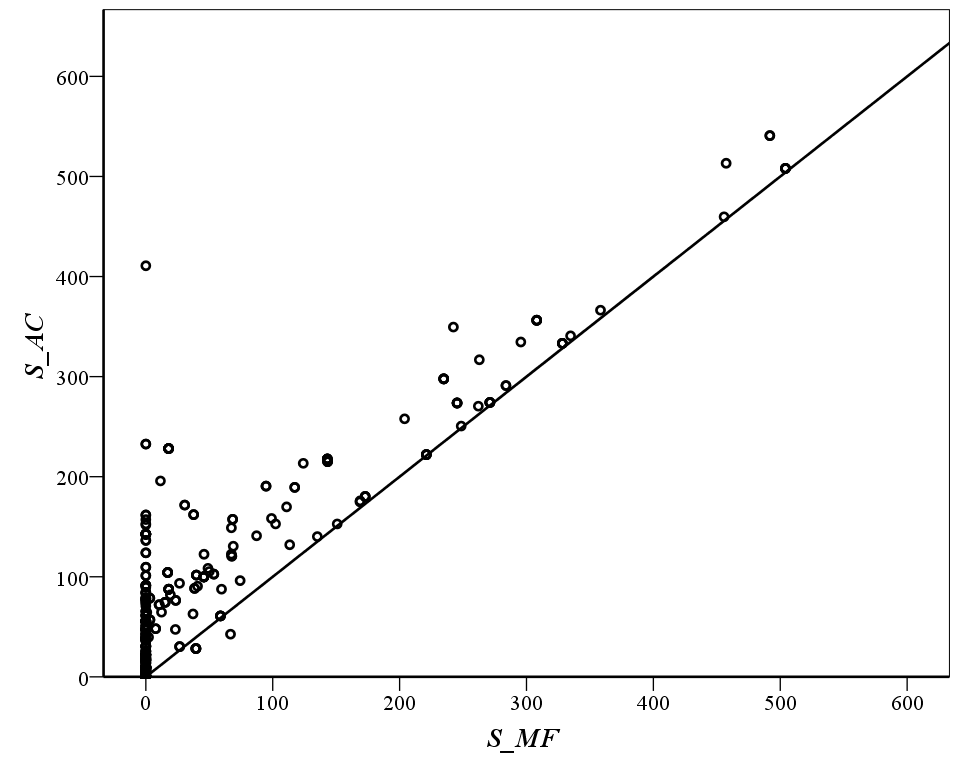}\label{fig:sl118}}
\subfigure[MEDIUM
system]{\includegraphics[width=0.49\linewidth]{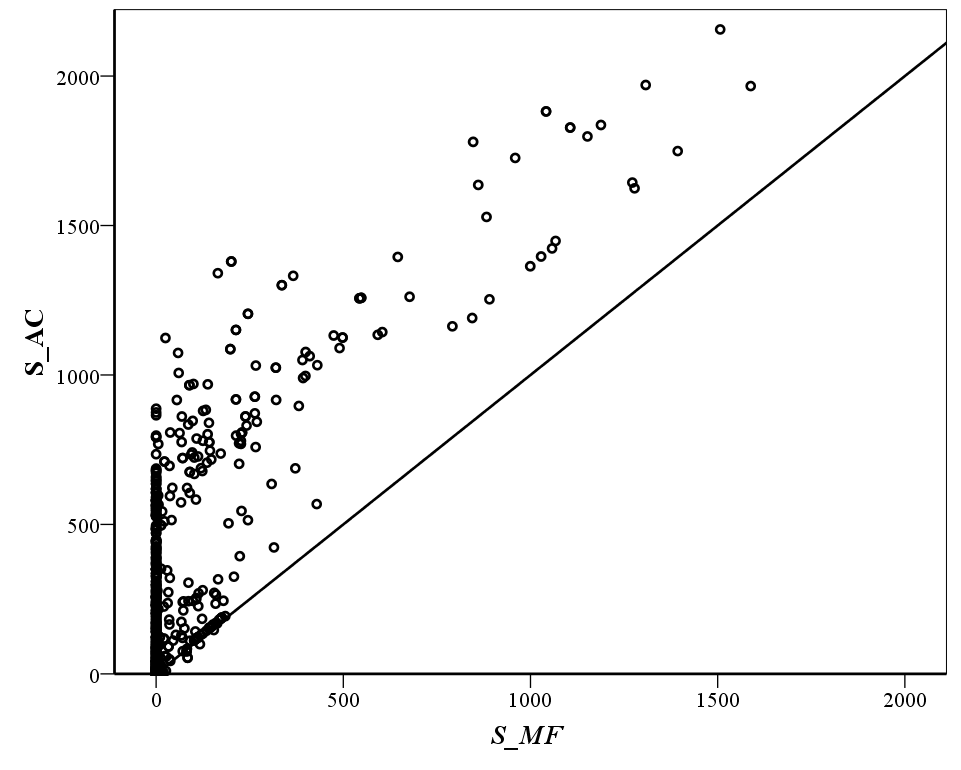}\label{fig:sl2383}}
%\subfigure[LARGE
%system]{\includegraphics[width=0.32\linewidth]{shedload9241.png}\label{fig:sl9241}}
\caption{Relation (scattering plot) between AC-OLS-based load
shedding $S_{AC}(s)$ (vertical axis) and maximum-flow-based load
shedding $S_{MF}(s)$ (horizontal axis). A diagonal line $y=x$ is
added for clarity.}\label{fig:shedload}
\end{figure*}

Figure~\ref{fig:shedload} shows the relation between $S_{MF}(s)$ and
$S_{AC}(s)$ for islands met in optimal partitions of SMALL and
MEDIUM powers systems. It can be seen from the figure that
$S_{MF}(\cdot)$ correlates well\footnote{Correlation is equal to
0.89 for SMALL system and 0.78 for MEDIUM system.} with $S_{AC}(s)$
but systematically underestimates the latter. Also, there are
numerous situations when $S_{MF}(s)=0$ while $S_{AC}(s)>0$.

Accuracy of $S_{AC}(\cdot)$ prediction can be improved sufficiently
by considering a three-variate non-linear regression
\begin{equation}\label{eq:S_AC_regression}
\bar{S}_{AC}(s) = \max\left[S_{MF}(s)-a GR_{MF}(s)+bw(s), 0\right],
\end{equation}
where $a$ and $b$ are regression parameters (grid-dependent, in
general), $w(s)$ is the volume of island $s$, while $GR_{MF}(s)$ is
the \emph{generation reserve}. The latter is calculated as
$GR_{MF}(s):=\sum_{i\in s} \left(G_i-g_i^*(s)\right)$, where $G_i$
is the maximum real power output of generator bus $i\in s$ and
$g_i^*(s)$ is its real power output in the solution of MF-OLS
problem for island $s$.

\begin{figure*}
\subfigure[SMALL
system]{\includegraphics[width=0.49\linewidth]{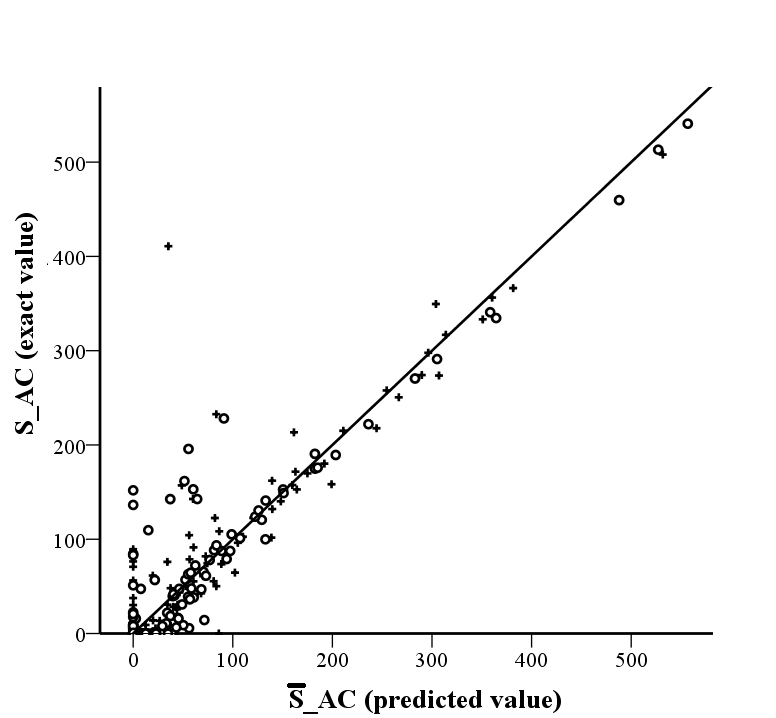}\label{fig:sl118_2}}
\subfigure[MEDIUM
system]{\includegraphics[width=0.49\linewidth]{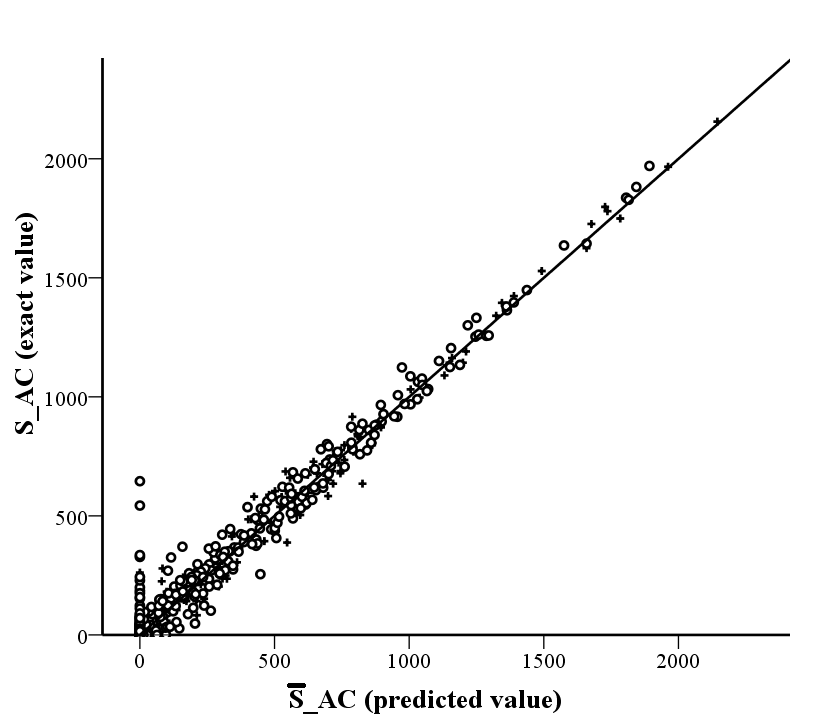}\label{fig:sl2383_2}}
%\subfigure[LARGE
%system]{\includegraphics[width=0.32\linewidth]{shedload9241.png}\label{fig:sl9241}}
\caption{Scattering plot for $S_{AC}(s)$ (vertical axis) and
 $\bar{S}_{AC}(s)$ (horizontal
axis). Circles denote points of the training set, while pluses go
for the testing set.}\label{fig:shedload_2}
\end{figure*}

For comparison, scattering plots (predicted value $\bar{S}_{AC}(s)$
\emph{vs} exact value $S_{AC}(\cdot)$) of regression
(\ref{eq:S_AC_regression}) under the best parameter values ($a=0.43,
b=0.02$ for SMALL system and $a=0.97, b=0.04$ for MEDIUM system) are
shown in Figures \ref{fig:sl118_2} and \ref{fig:sl2383_2}
respectively. Correlation is 0.92 with mean absolute error (MAE)
13.9 MW on the training set and is 0.94 with MAE 14.1 on the testing
set\footnote{Training and testing sets were obtained with halfway
random sampling.} for SMALL system. For MEDIUM system correlation is
equal to 0.99 with MAE 39.2 MW on the training set and is equal to
0.98 with MAE 39.5 MW on the testing set.\footnote{Distinct to SMALL
system, the definition of MEDIUM system includes realistic real
power flow constraints for transmission lines, which results in
better prediction accuracy.}

To take advantage of the more accurate regression
(\ref{eq:S_AC_regression}) only minor change in MIQP setting
(\ref{eq:MIQP}) is needed. It is enough to introduce the new
non-negative continuous variables $\sigma_k\ge 0$, $k=1,...,K$,
which satisfy inequality constraints
$$\sigma_k\ge \sum_{i=1}^n x_{ik}\left[l_i-a(G_i-g_i)+bw_i\right],k=1,...,K$$
and define the new cost function
\begin{equation}\label{eq:MIQPcost2}
\tr X^T L(A)X + \sum_{k=1}^K \sigma_k.
\end{equation}

The revised MIQP still has linear constraints and the convex cost
function. Therefore, its complexity does not increase, while its
solution accurately predicts the optimal value of the cost function
$$F(\pi)=\alpha_C C(\pi)+\alpha_D D(\pi)+\alpha_S S_{AC}(\pi)$$
in the considered test systems and is expected to meet
application-specific accuracy requirements.

The final aspect of algorithm accuracy is that MIQP (\ref{eq:MIQP}) is solved only
approximately by ISC algorithm due to high problem dimensionality, and the
approximation error of ISC algorithm should be estimated by comparing the ISC
solution to the exact solution of MIQP (\ref{eq:MIQP}). Due to computational
intractability of problem (\ref{eq:MIQP}) for large grids, such a direct error
evaluation is impossible for MEDIUM and LARGE system, while for SMALL system the
mean relative error of ISC algorithm for 100 cases is just $2\%$ (the maximum error
is $9\%$). Therefore, approximation error is negligible.

%Spectral clustering is applied first to produce the aggregated
%low-dimensional graph, which is partitioned with the branch-and-cut
%scheme of CPLEX solver.
Let us continue with algorithm benchmarking. The proposed ISC
algorithm is based on CSC and HSC, the most efficient
spectral-clustering-based algorithms for OCI. It is natural to use
them as a benchmark in performance evaluation of ISC for large grids
where exact solution is intractable. Unfortunately, $K$-partition
$\chi_{K}(|P|~|\mathbf{w})$ computed by HSC algorithm is often
imbalanced, and so does $K$-partition
$\chi_{K}(\tilde{\Phi}|\mathbf{w})$, which we used to elicit
coherent generator groups for CSC.

\begin{figure*}
\subfigure[SMALL
system]{\includegraphics[width=0.5\linewidth]{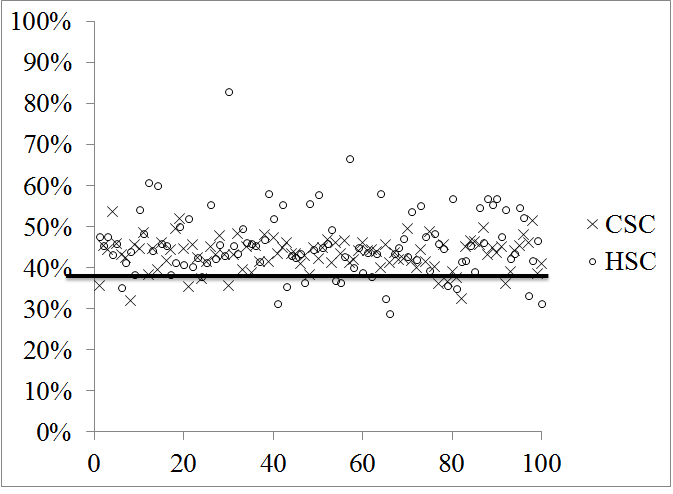}\label{fig:imbalance118}}
\subfigure[MEDIUM
system]{\includegraphics[width=0.5\linewidth]{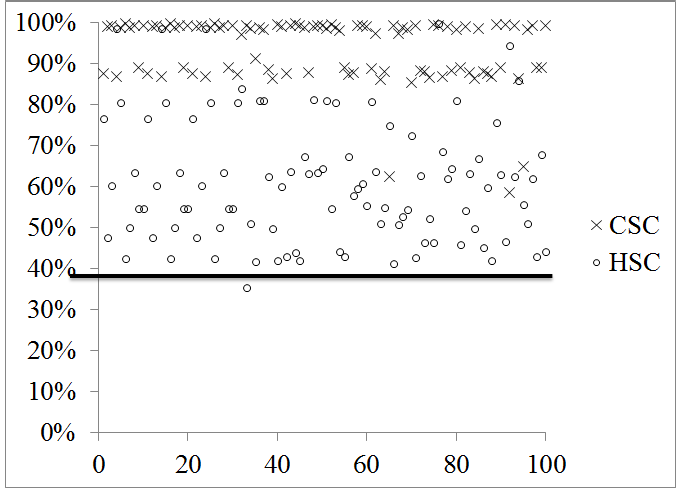}\label{fig:imbalance2383}}
\caption{Relative volume (per cent) of the largest island for
$K$-partitions generated by CSC and HSC. Case numbers are depicted
along the horizontal axis.}\label{fig:imbalance}
\end{figure*}

As a consequence, both CSC and HSC applied directly to our data in
most cases return imbalanced partitions being highly impractical for
controlled islanding. Figure \ref{fig:imbalance} shows the volume of
the largest island (relative to the total grid volume) for SMALL and
MEDIUM test systems. The horizontal line shows the maximum volume
constraint. HSC partitions are shown with circles, while crosses
stand for CSC partitions. It follows from the figure that CSC
returns a balanced partition only in 18 cases of 100 for SMALL
system, while HSC results in a balanced partition in 12 cases. For
MEDIUM system HSC returns a balanced condition just once, while CSC
partitions are always imbalanced.

Win rates of seven strategies of ISC algorithm for three test
systems are presented in Table \ref{tab:winrate}. The winning
strategy minimizes the cost function $F_{AC}(\cdot)$ with shed load
evaluated from AC-OLS problem. Strategies \ref{strategy_fixedPA} and
\ref{strategy_fixedP} appear the most successful for all three test
systems, while Strategies \ref{strategy_minP} and \ref{strategy_CSC}
only occasionally win (These strategies represent the closest
balanced analog of HSC and CSC algorithms consequently, so they can
be, to some extent, used as a baseline to compare ISC to the
competing algorithms.)

\begin{table}[!ht]
\caption{Win rate (per cent) of different strategies of ISC
algorithms. Win rates of Strategies \ref{strategy_minP} and
\ref{strategy_CSC} (separated by horizontal lines) can be used as a
baseline for comparison of ISC with, respectively, HSC and CSC
algorithms.} \label{tab:winrate}
\begin{tabular}{llrrr}
\hline\noalign{\smallskip}
   & & \multicolumn{3}{c}{Test power system} \\
  \multicolumn{2}{l}{Strategy} & SMALL & MEDIUM & LARGE \\
\noalign{\smallskip}\hline\noalign{\smallskip}
  \ref{strategy_fixedA} & Fixed-granularity strategy for  $A$    & 21 & 9  & 1  \\
  \ref{strategy_fixedP} & Fixed-granularity strategy for $|P|$   & 16 & 33 & 32 \\
  \ref{strategy_minA}   & Minimum-granularity strategy for $A$   & 20 & 10 & 4  \\
\hline
  \ref{strategy_minP}   & Minimum-granularity strategy for $|P|$ & 3  & 4  & 14 \\
  \ref{strategy_CSC}    & Minimum-granularity refined CSC        & 0  & 2  & 3  \\
\hline
  \ref{strategy_fixedPA}& Sequential strategy for $|P|$ and $A$  & 29 & 37 & 38 \\
  \ref{strategy_meet}   & Crossing CSC and HSC partitions        & 11 & 5  & 8 \\
\noalign{\smallskip}\hline
\end{tabular}
\end{table}

\begin{figure*}
\subfigure[SMALL
system]{\includegraphics[width=0.5\linewidth]{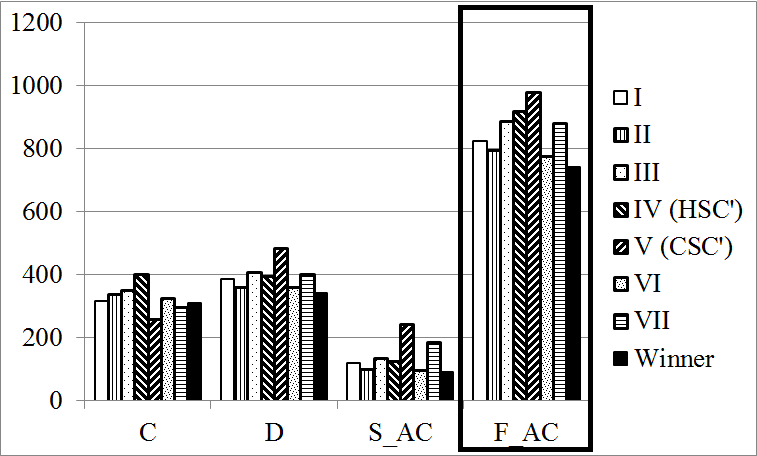}\label{fig:118s}}
\subfigure[MEDIUM
system]{\includegraphics[width=0.5\linewidth]{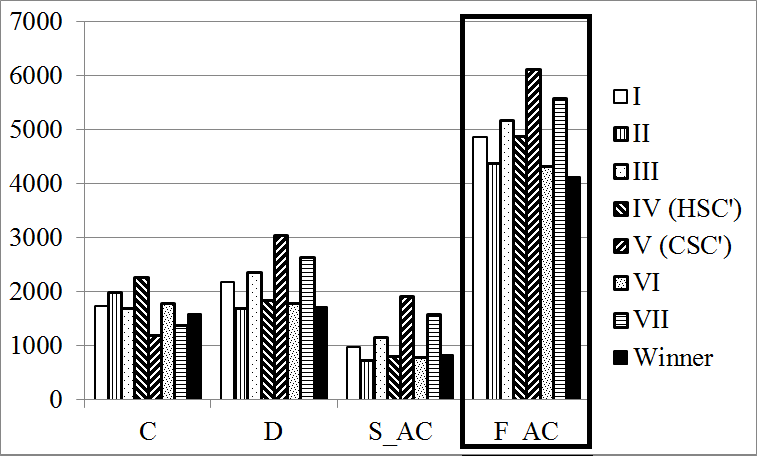}\label{fig:2383s}}
\caption{Average performance metrics for Strategies
\ref{strategy_fixedA}-\ref{strategy_meet} of ISC algorithm. Total
cost $F_{AC}(\cdot)$ is framed. Generator coherence $C(\cdot)$, flow
disruption $D(\cdot)$, and shed load $S_{AC}(\cdot)$ are also
shown.}\label{fig:performance}
\end{figure*}

The average partitioning quality of ISC algorithm for SMALL test
system is depicted in Figure \ref{fig:118s}. The value of cost
function is enclosed in a frame, all its components are also
presented for all seven strategies and for the winning strategy.
Although Strategy \ref{strategy_fixedP} (fixed granularity for
matrix $|P|$) does not account directly for generator coherency, it
suggests competitive coherency cost and wins in disruption and load
shedding. At the same time, Strategy \ref{strategy_fixedPA}
(sequential partitioning) has the same average cost due to a little
bit lower generator coherency under a slightly higher disruption and
shed load. These two strategies remain the leaders for all three
test systems.\footnote{This result is valid for the fixed weights of
performance metrics: $\alpha_C=\alpha_D=\alpha_S=1$. Having the
weights changed, the leading strategies can also change (see the
analysis in Section \ref{sec:flex}).}

A yet another conclusion is that considering excessive number of
eigenvectors (which is the main idea of ISC algorithm) is essential
to achieve good performance, since Strategy \ref{strategy_minP},
which differs from Strategy \ref{strategy_fixedP} only in the number
of eigenvectors considered, is inferior. The same is true for
Strategy \ref{strategy_CSC}, which is based on CSC and also operates
with the limited number of eigenvectors. Performance metrics of ISC
strategies for MEDIUM system are presented in Figure
\ref{fig:2383s}.

The similar analysis for LARGE system is hindered by the fact that
for this system some strategies often suggest imbalanced partitions
(Strategies \ref{strategy_fixedA}, \ref{strategy_CSC}, and
\ref{strategy_meet} show the worst success rate: 24\% , 32\%, and
34\% respectively) even when $4K$ eigenvectors are calculated. Since
island volume constraints are mandatory, ISC algorithm neglects
imbalanced partitions suggested by some strategies when selecting
the final solution, and there is no common base for the comparison
of these strategies.

Therefore, maintaining several different partitioning strategies is
important to always obtain a good islanding solution in large power
systems.

\subsection{Flexibility}\label{sec:flex}

\begin{figure*}
\subfigure[]{\includegraphics[width=0.5\linewidth]{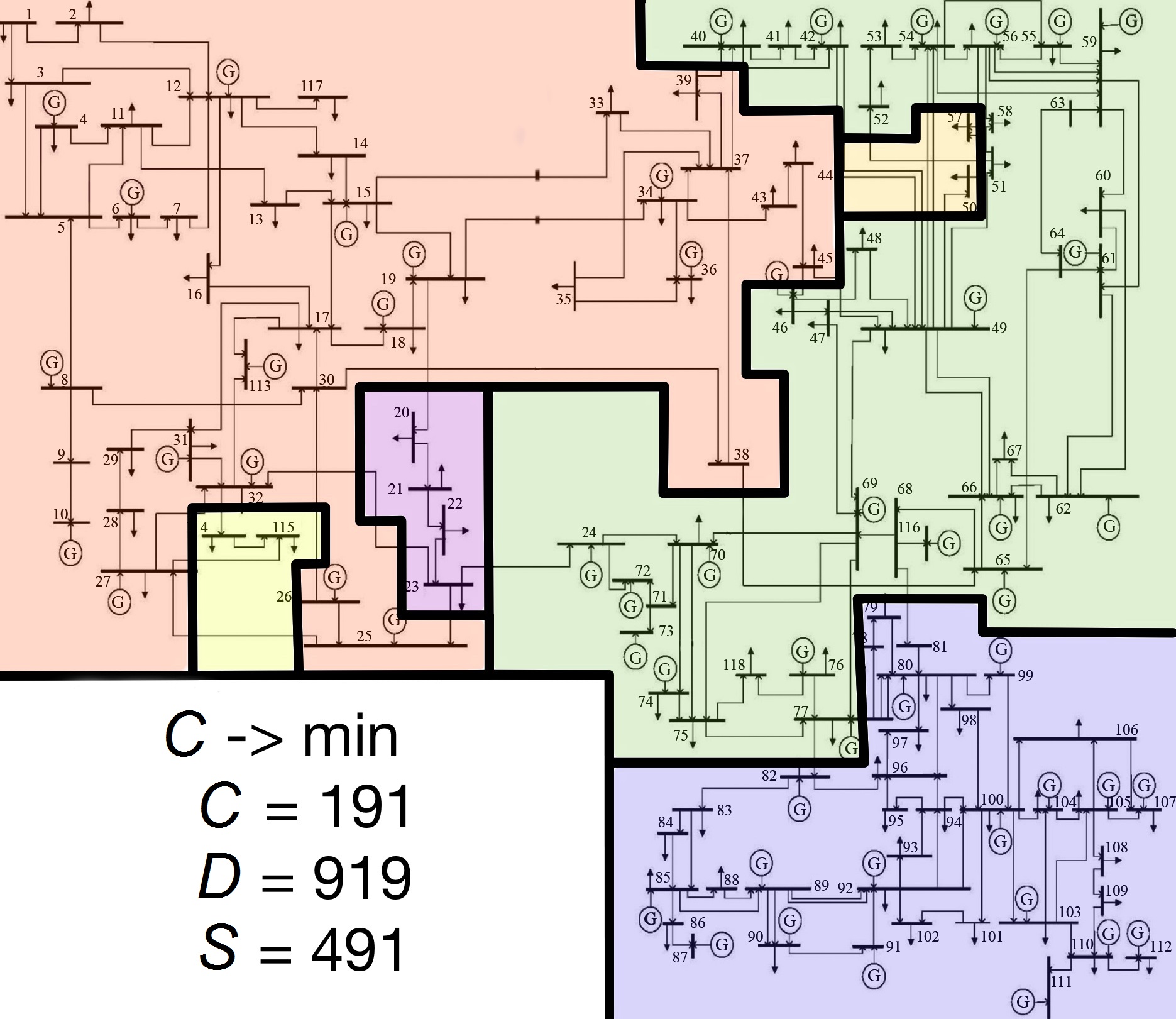}\label{fig:costvarC}}
\subfigure[]{\includegraphics[width=0.5\linewidth]{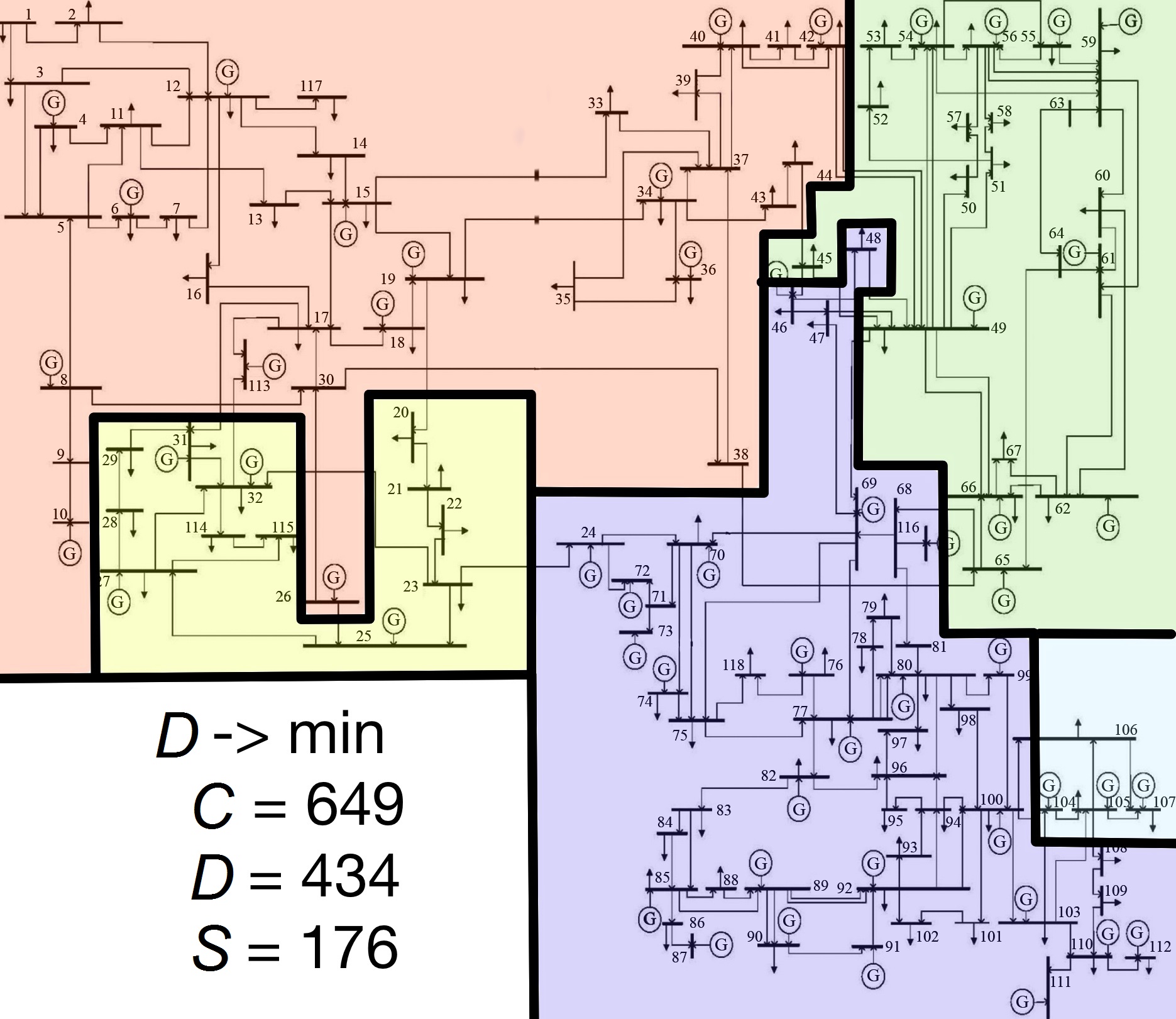}\label{fig:costvarP}}\\
\begin{center}
\subfigure[]{\includegraphics[width=0.5\linewidth]{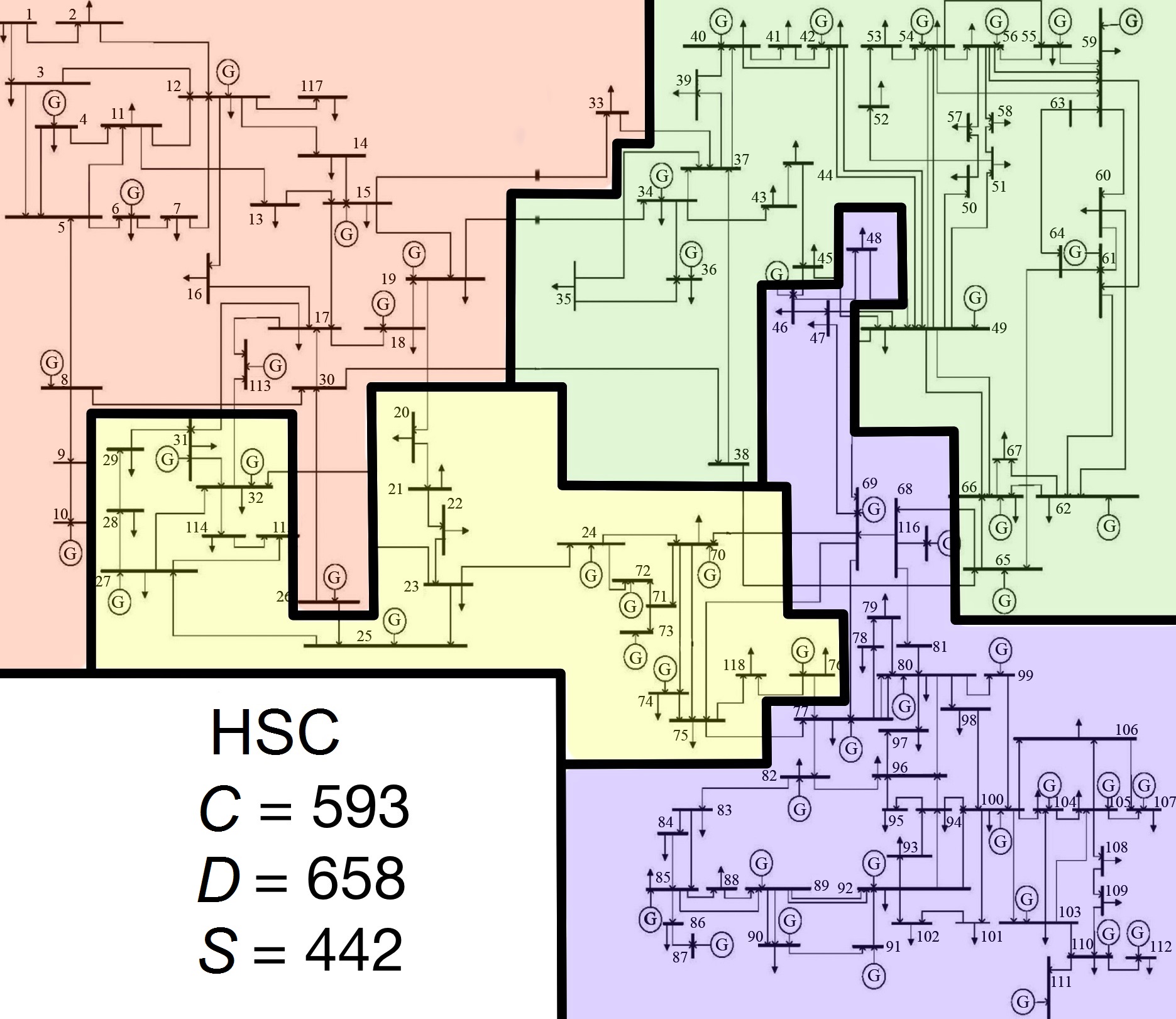}\label{fig:costvarbase}}\\
\end{center}
\subfigure[]{\includegraphics[width=0.5\linewidth]{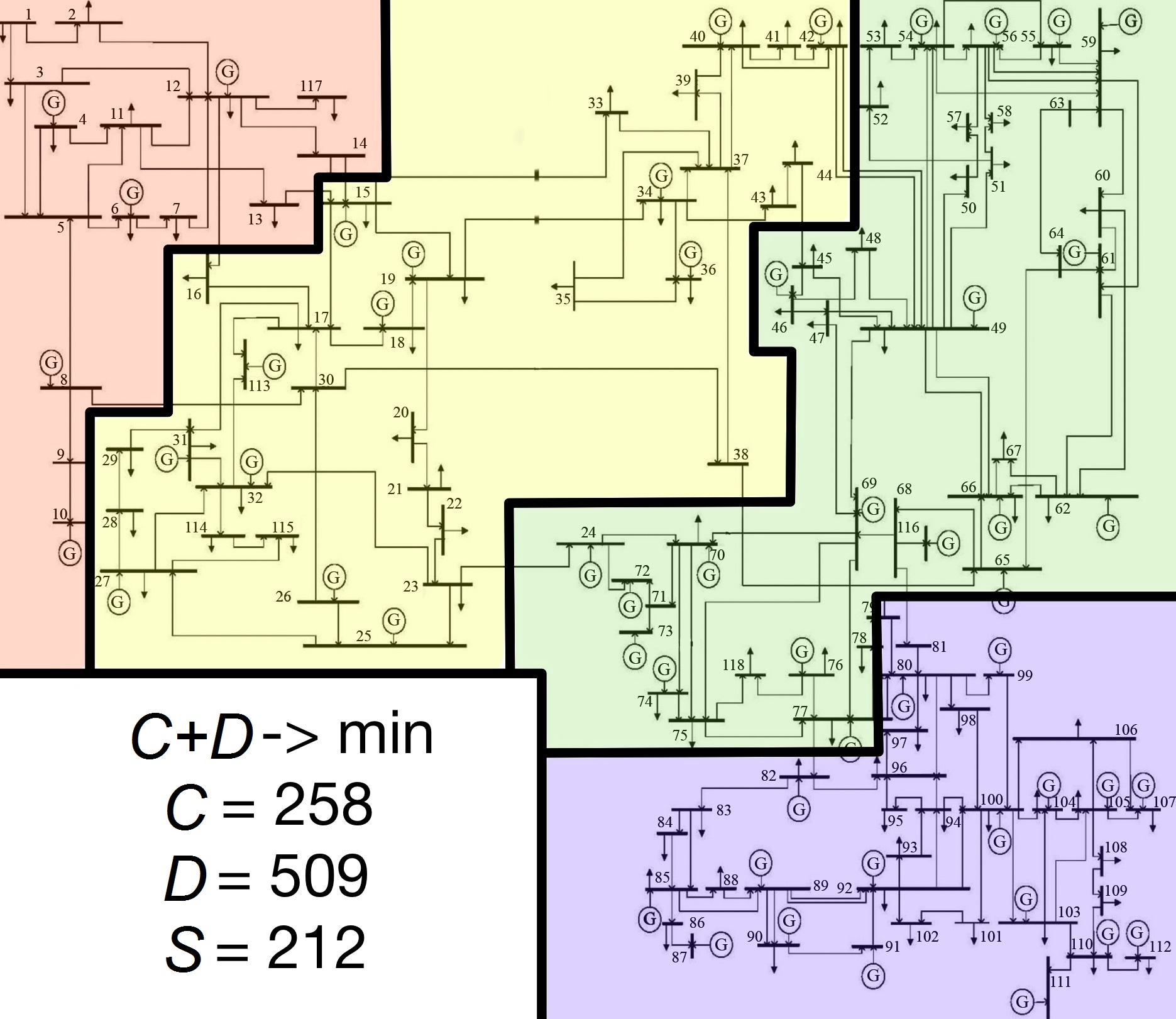}\label{fig:costvarCP}}
\subfigure[]{\includegraphics[width=0.5\linewidth]{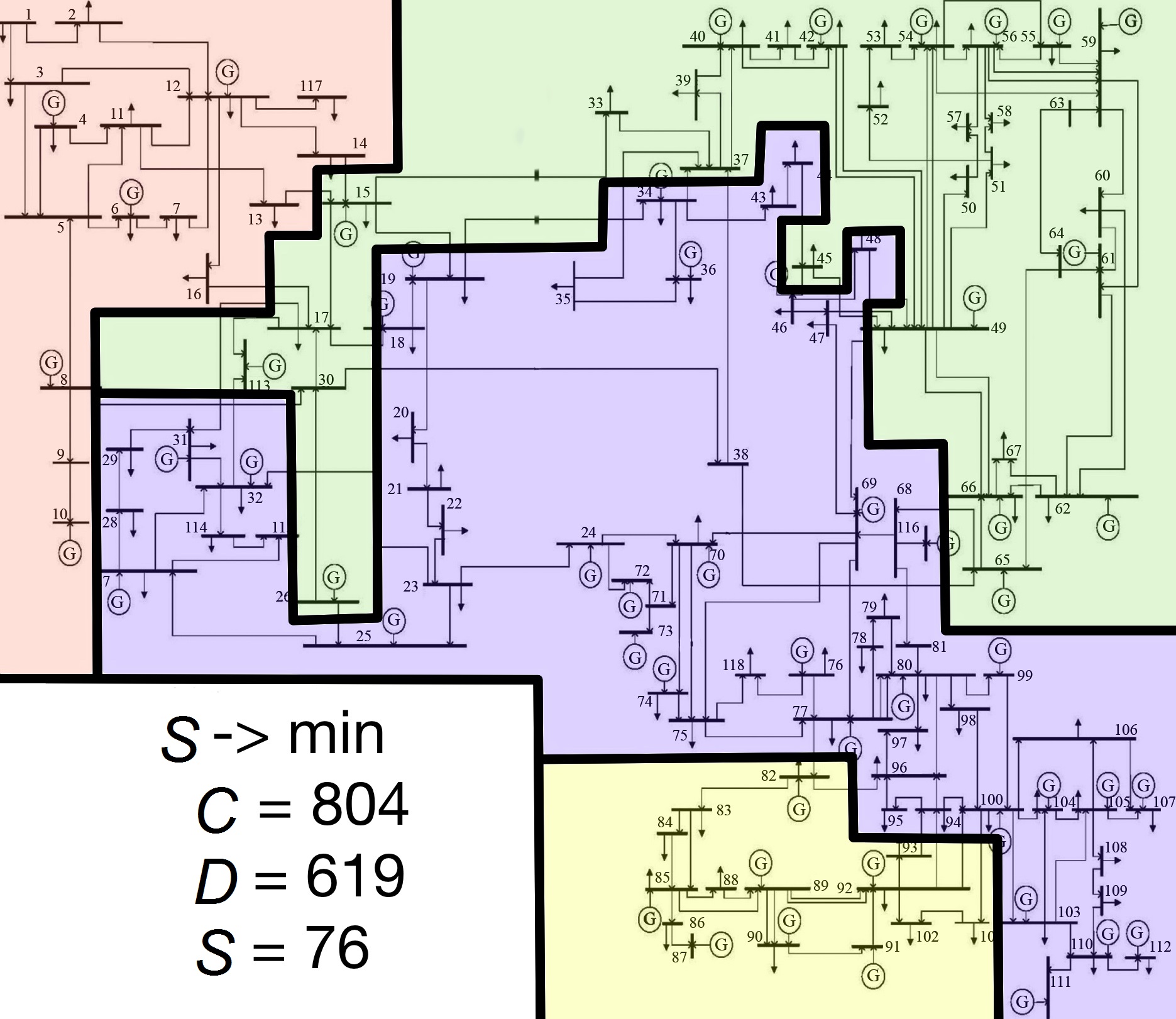}\label{fig:costvarS}}
\caption{Optimal islanding schemes for different cost functions} \label{fig:costvar}
\end{figure*}

In many existing OCI algorithms
\cite{li2005strategic,xu2010SimplifyMETIS,ding2013CoherencyDisruptionConstrainedSpectralBisection,sanchez2014AdmittanceDisruptionHierSpec,quiros2015constrained}
an optimization criterion cannot be tuned to assign different
importance to different performance metrics of controlled islanding.
Distinct to them, in ISC algorithm the weights of cost function
components (dynamic coupling, disruption, and shed load) can be
chosen arbitrary leading to different resulting partitions. For
example, in Figure \ref{fig:costvar} optimal islanding schemes for
SMALL system under several extremal settings are presented. In
Figure \ref{fig:costvarC} weights in cost function (\ref{eq:F_def})
were chosen to minimize dynamic coupling $C(\cdot)$ while ignoring
$D(\cdot)$ and $S(\cdot)$. On the contrary, in Figure
\ref{fig:costvarP} disruption $D(\cdot)$ is minimized with no
attention paid to $C(\cdot)$ and $S(\cdot)$. The partition that
balances generators' dynamic coupling $C(\cdot)$ and disruption
$D(\cdot)$ is presented in Figure \ref{fig:costvarCP}, and the one
taking care only for the shed load $S(\pi)$ is presented in Figure
\ref{fig:costvarS}.

For comparison, the best alternative partition (calculated by HSC
algorithm) and its performance metrics are presented in
Figure~\ref{fig:costvarbase}. Therefore, every single performance
metric of HSC partition or a combination of metrics can be improved
by ISC algorithm, and the optimization criterion can be flexibly
tuned to adjust relative importance of different metrics.

% more analysis on specific partitions

\subsection{Computational Complexity}

%\textbf{What is needed to be calculated: $r_1 K A, r_2K |P|, r_3K
%\tilde{\Phi}$. a series of bisections for Sequential strategy. The first is the hardest, max dimension, max density.
%}

Compared to the competitive spectral partitioning algorithms (CSC
and HSC), the proposed ISC algorithm requires additional
calculations: in the first step several overdetailed spectral
partitions are calculated (strategies I-VII are introduced in
Section \ref{sec:Spectral}), and in the second step every detailed
partition $\pi^\text{I}, ..., \pi^\text{VII}$ is merged into a
$K$-partition by CPLEX optimization
routines.\footnote{\texttt{cplexmiqp} routine of CPLEX 12.6.2.0 was
used to solve MIQP; tests were run on Intel Core i-5 3337U CPU 1.8
GHz.} Therefore, total computation time of ICS algorithm is equal to
the maximum processing time for partitioning strategies I, ..., VII
(although the first suggestion is given as soon as the fastest
strategy completes).

%
% ADD Here
% Some eigenvectors are demanded by several Strategies.
% So, they can be computed once and reused. Figure X shows interdependencies
% of parallel computations for ISC algorithm. Box sizes are nominal and are not related directly to operation duration.
%

\begin{figure*}
  \includegraphics[width=1\textwidth]{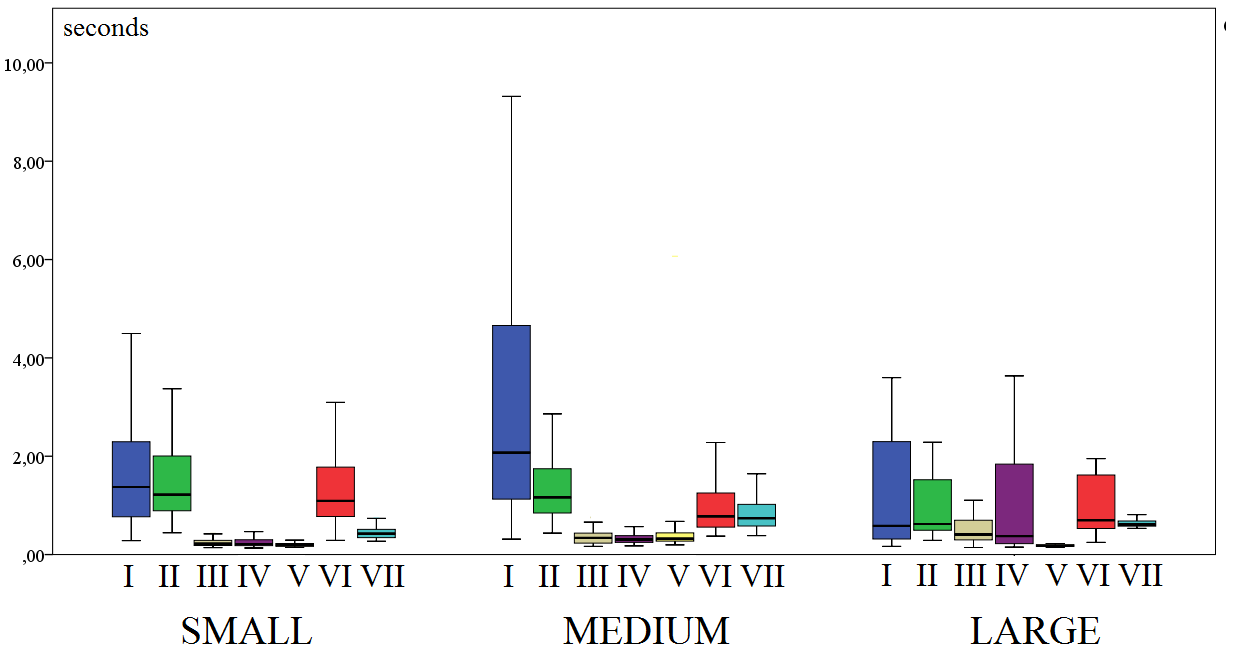}
\caption{Time complexity of MIQP (the second step o ISC algorithm)
for 3 test systems and 7 partitioning strategies. The median time
and quartile boundaries are presented.} \label{fig:step2}
\end{figure*}

The detailed statistics of MIQP solution time is shown in Figure
\ref{fig:step2}. Solution time is lower for minimum-granularity
strategies III-V due to lower dimensionality of the problem. On the
other hand, solution time depends insignificantly on the test system
as the order of an aggregated graph does not depend on the size of
an initial grid but only on the number of connected islands $n'$ in
the detailed partition, which linearly depends on the number $K$ of
islands requested (remember that $K=4$ in all experiments).

%\footnote{For
%MEDIUM system computational time is higher, probably, due to its
%higher connectivity.}

%Higher volatility of computation time is observed for Strategies I,
%II, VI. As Figure \ref{fig:performance} shows, These strategies are
%the most efficient ones.

\begin{figure*}
  \includegraphics[width=0.7\textwidth]{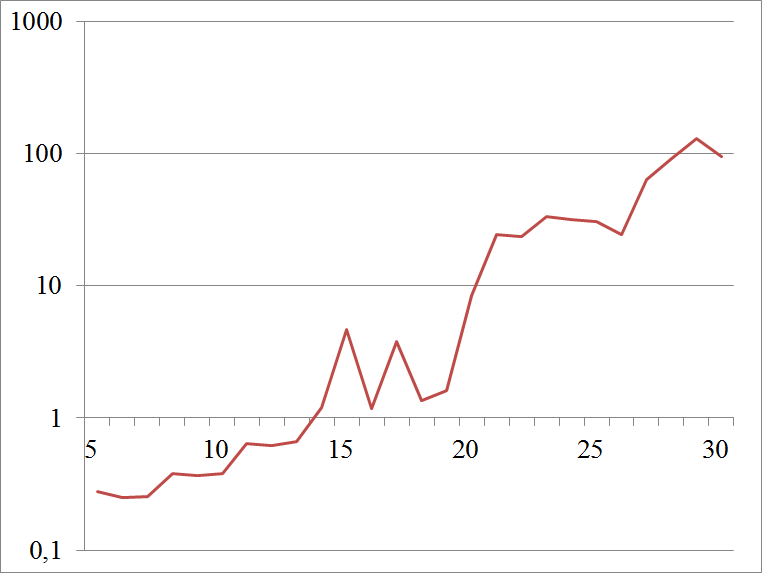}
\caption{Time complexity of MIQP vs order of aggregated graph being
partitioned (logarithmic scale).} \label{fig:MIQPdimension}
\end{figure*}

Figure \ref{fig:MIQPdimension} presents the average MIQP computation
time as a function of dimensionality of aggregated graph $n'$ (SMALL
system is used for this test). Computation time grows exponentially,
therefore, for $n'>K_{max}\approx 20$ network dimension should be
reduced by the greedy algorithm before running MIQP.

\begin{figure*}
\subfigure[SMALL
system]{\includegraphics[height=0.35\linewidth]{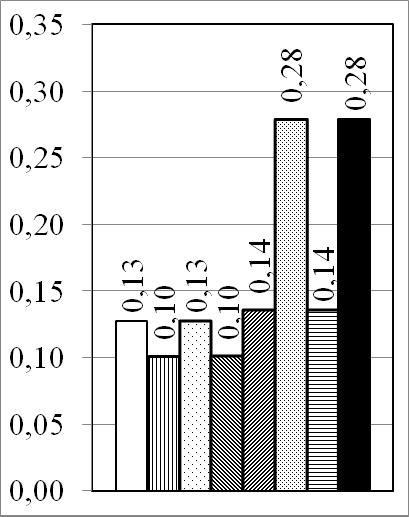}\label{fig:compISCSMALL}}
\subfigure[MEDIUM
system]{\includegraphics[height=0.35\linewidth]{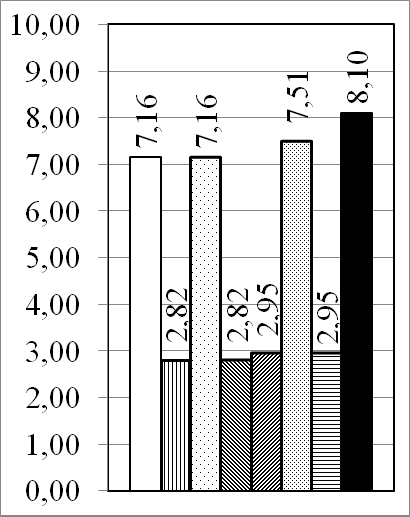}\label{fig:compISCMEDIUM}}
\subfigure[LARGE
system]{\includegraphics[height=0.35\linewidth]{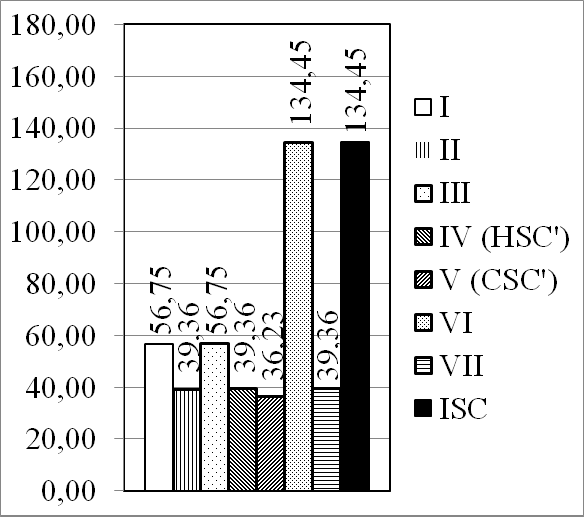}\label{fig:compISCLARGE}}
\caption{Average time of eigenvectors' calculation (the first step
of ISC) for 7 parallel strategies of ISC algorithm. Eigenvectors
were evaluated with \texttt{eigs} routine of Matlab R2014a run
on Intel Core i-5 3337U CPU (1.8 GHz) laptop.} \label{fig:compISC}
\end{figure*}

Figure \ref{fig:compISC} illustrates average time complexity of the
first step of ISC algorithm for all seven parallel strategies.
Comparing to Figure \ref{fig:step2}, we conclude that the second
step of ISC algorithm is fast enough compared to the first step in
MEDIUM and LARGE systems, and its computation time is significant
only in SMALL system.

The second step of ISC algorithm can be fostered in two ways. The
first issue is high variability of computation time (A typical
problem of branch-and-bound procedures.) Sometimes the exact
solution takes much longer than usually, and it may be critical for
online OCI applications. This problem is solved by imposing a sharp
time limit on MIQP calculations. Fortunately, in ISC algorithm the
branch-and-bound procedure always starts with a feasible record
solution produced by the greedy algorithm. The greedy algorithm is
fast enough (less than 0.03 sec. on average) and often provides
good solutions (the average cost improvement of the exact algorithm
over the greedy one varies from 1.5\% for LARGE system to 7\% for
SMALL system).

The second approach to foster MIQP calculations assumes increasing
the number of processing units. Branch-and-bound procedures are
easily parallelized, and doubling the number of processors almost
doubles the performance.

Figure \ref{fig:compISC} shows that evaluation of eigenvectors of
the normalized Laplace matrix is the most computationally expensive
operation of spectral clustering algorithms (at least, for
realistically sized grids). Its complexity depends on the number of
eigenvectors requested, on matrix dimension, and on its sparsity
ratio.

In Strategies I and III $r_1K$ eigenvectors of matrix
$L_{sym}(A|\mathbf{w})=L_{sym}(\alpha_C\Phi+\alpha_D|P|~|\mathbf{w})$
are calculated. Strategies II and IV require calculation of $r_2K$
eigenvectors of $L_{sym}(|P|~|\mathbf{w})$. The latter matrix is
sparser and, therefore, the first step of Strategies I and III is
computationally more expensive than that of Strategies II and IV.

Strategy V (CSC-based algorithm) is slightly more computationally
expensive than Strategy IV (HSC-based) for SMALL and MEDIUM systems
as it involves calculation of the eigenproblem for the Laplacian of
denser matrix $\tilde{\Phi}$ of generators' dynamic coupling. At the
same time, for LARGE system Strategy V is faster than Strategy IV
due to the smaller size of matrix $\tilde{\Phi}$.

Strategy VI involves the nested spectral bisection of $K$-partition
and appears the most time-consuming strategy for all three test
systems. On the other side, it is the most successful one (according
to Table \ref{tab:winrate}) and is very reliable in generating
admissible (balanced) partitions. Comparing Figures
\ref{fig:compISC} and \ref{fig:performance} we see that Strategy II
is approximately three times faster than Strategy IV and results in
only minor loss of partitioning quality. Therefore, one may consider
abandoning Strategy VI if it fails to satisfy application-specific
time limits.

Strategy VII takes partitions $\pi_\text{IV}$ and $\pi_\text{V}$ as
an input, so its computation time is the maximum of those for
Strategies IV and V (computation time of the meet operation can be
neglected).

%Distinct to ISC algorithm, only spectral $K$-partitions are
%evaluated and $K$ eigenvectors are calculated in CSC and HSC
%algorithms.\footnote{Nevertheless, to obtain balanced partitions
%that always satisfy island volume constraints more granulated
%partitions can be requested, which requires more eigenvectors.}

%Practical computation time of CSC, HSC, and ISC algorithms in three
%test grids is presented in Figure \ref{fig:complex}. As before,
%$K=4$ is the requested maximum number of islands, $r_1=r_2=r_4=4$,
%so, detailed partitions have 16 islands. All algorithms where
%implemented in Matlab R2014a and run on Intel Core i-5 3337U CPU
%(1.8 GHz).

%\begin{figure*}
%  \includegraphics[width=0.9\textwidth]{complexity.png}
%\caption{Time complexity of algorithms for three grids}
%\label{fig:complex}       % Give a unique label
%\end{figure*}

For SMALL system eigenvector calculation is very fast, and, hence,
ISC algorithm is slower than CSC \cite{quiros2015constrained} and
HSC \cite{sanchez2014AdmittanceDisruptionHierSpec} algorithms
(mostly, due to the additional second step, which takes about 2 sec.
irrespective of the system dimension). For MEDIUM system the speed
is comparable (15 seconds is reported for HSC
\cite{sanchez2014AdmittanceDisruptionHierSpec}, while CSC algorithm
was tested only for the reduced Great Britain 815-bus system in
\cite{sanchez2014AdmittanceDisruptionHierSpec}). These algorithms
cannot be directly compared for LARGE system, since CSC and HSC fail
to provide eligibly balanced partitions. If we relax the maximum
island volume constraint, we can expect ISC to be at least three
times slower than CSC and HSC, mostly due to the slow Strategy VI
(which, as noted above, can be abandoned in case of time deficit).

It should be noted that for LARGE system (and, perhaps, for MEDIUM system) neither
of existing spectral-clustering-based algorithms is fast enough for online
islanding, when run on a PC. Eigenvectors' calculation takes most time, and,
fortunately, it can be dramatically fostered by using parallel vector calculus
capabilities of graphical processing units (GPU)
\cite{lessig2007eigenvalue,bell2008efficient}. Hence, fast implementation of
eigenvector calculation for spectral clustering OCI algorithms in large power
systems is no more than a programming issue. Potentially, the performance
improvement can also be obtained by replacing spectral clustering with some
algorithms of $NCut$ minimization based on semidefinite programming
\cite{fan2012multi,ames2014guaranteed}. A yet another promising approach is to use
deep learning techniques to replace the optimization problem solution with fast
calculation by a multilayered artificial neural network. In this case the algorithm
proposed in this article calculates an extensive data set used to learn the neural
network. Comparison of these approaches requires additional work.

\section{Conclusion}\label{sec:Conclusion}
An improved version of spectral-clustering-based OCI algorithms
\cite{sanchez2014AdmittanceDisruptionHierSpec,quiros2015constrained} is proposed in
this article, which accounts for generator coherency, power flow disruption, shed
load, and island size, allowing to adjust flexibly relative importance of these
performance metrics and caring for the maximum island volume.

An optimal islanding scheme is sought in two steps. Firstly, several
alternative detailed grid partitions are determined that contain
more islands than required while balancing generator coherency and
power flow disruption. Secondly, some islands in these partitions
are merged to optimize the weighted sum of generator coherency,
power flow disruption, and shed load while fulfilling the maximum
island volume constraint. The second step reduces to MILP whose
dimension does not depend on the dimension of the original grid and
is small enough to use exact algorithms.

A series of experiments for three standard test grids was performed.
The results show that, compared to competing algorithms, HSC
\cite{sanchez2014AdmittanceDisruptionHierSpec} and CSC
\cite{quiros2015constrained}, the proposed improved spectral
clustering algorithm results in substantial improvement of both the
overall composite criterion of partition quality (more than twice
for some grids) and of each its component. The algorithm is
computationally efficient. Potentially, it can be used to partition
large power systems in real time.

The proposed algorithm can also be helpful in applications other
than energy systems to partition high-dimensional directed graphs
with asymmetric matrix of vertex weights.

The multi-objective criterion studied is a weighted sum of various performance
metrics (generators' dynamic coupling, power flow disruption, ECI, excess demand,
shed load, and, probably, others), but the choice of performance metrics included
into the optimization criterion and assignment of their relative weights is an open
question in general. Further research is needed to establish an empirically grounded
relation between transient stability of an island being created (e.g., with
simulation techniques from \cite{ames2014guaranteed,song2016dynamic}) and
computationally efficient performance metrics (some of which were mentioned above).

\begin{acknowledgements}
The first author would like to thank support from Russian Foundation
for Basic Research (project 16-37-60102).
\end{acknowledgements}

\bibliographystyle{spmpsci}      % mathematics and physical sciences
\bibliography{GoubkoGinz}   % name your BibTeX data base

\end{document}